\documentclass[12pt]{article}
\usepackage[latin1]{inputenc}
\usepackage{tgtermes} 
\usepackage[T1]{fontenc}
\usepackage[english]{babel}
\usepackage{amssymb,amsfonts,amsmath}
\usepackage{authblk}
\usepackage{cite}
\usepackage{float}

\usepackage{booktabs}
             
\usepackage{geometry} 
\usepackage{multirow} 
\usepackage[hyperindex,colorlinks=false]{hyperref} 
\usepackage{color} 
\usepackage{graphicx} 
\usepackage{subfig}
\usepackage{tabularx}
\usepackage{tikz}
\usetikzlibrary{shapes}
\usepackage{xcolor}
\usepackage[normalem]{ulem} 
\DeclareMathOperator*{\argmin}{arg\,min}

\title{A Deep Collocation Method for the Bending Analysis of Kirchhoff Plate}

\author[a]{Hongwei Guo}
\author[b]{Timon Rabczuk}
\author[a,c,d]{Xiaoying Zhuang  \footnote{Corresponding authors:\\ \mbox{Xiaoying Zhuang}, \mbox{+49 511 762-19589}, \mbox{Xiaoying.Zhuang@gmail.com}}}
\affil[a]{Institute of Continuum Mechanics, Leibniz Universit\"{a}t Hannover, Appelstra{\ss}e 11, 30157~Hannover,~Germany}
\affil[b]{Institute of Structural Mechanics,~Bauhaus-Universit\"{a}t  
Weimar,~Marienstr.15~D-99423~Weimar,~Germany}
\affil[c]{Department of Geotechnical Engineering, Tongji University, Siping Road 1239, 200092~Shanghai,~P.R.China}
\affil[d]{Key Laboratory of Geotechnical and Underground Engineering of Ministry of Education, Tongji University, 200092~Shanghai, ~P.R.China}

\date{}
\begin{document}

\maketitle

\abstract{
\noindent In this paper, a deep collocation method (DCM) for thin plate bending problems is proposed. This method takes advantage of computational graphs and backpropagation algorithms involved in deep learning. Besides, the proposed DCM is based on a feedforward deep neural network (DNN) and differs from most previous applications of deep learning for mechanical problems. First, batches of randomly distributed collocation points are initially generated inside the domain and along the boundaries. A loss function is built with the aim that the governing partial differential equations (PDEs) of Kirchhoff plate bending problems, and the boundary/initial conditions are minimised at those collocation points. A combination of optimizers is adopted to in the backpropagation process to minimize the loss function so as to obtain the optimal hyperparameters. In Kirchhoff plate bending problems, the C1 continuity requirement poses significant difficulties in traditional mesh-based methods. This can be solved by the proposed DCM, which uses a deep neural network to approximate the continuous transversal deflection, and is proved to be suitable to the bending analysis of Kirchhoff plate of various geometries.}


\section{Introduction}

Thin plates are widely employed as basic structural components in engineering fields \cite{ventsel2001thin}, which combines light weight, efficient load-carrying capacity, economy with technological effectiveness. Their mechanical behaviours have long been studied by various methods such as finite element method \cite{bathe2006finite,hughes2012finite}, boundary element method \cite{katsikadelis2016boundary,brebbia2016boundary}, meshfree method \cite{liu2009meshfree}, isogeometric analysis \cite{nguyen2015isogeometric}, and numerical manifold method \cite{zheng2013numerical,guo2018linear,guo2019numerical}. The Kirchhoff bending problem is a classical fourth-order problem, its mechanical behaviour is described by fourth-order partial differential equation as it is pretty difficult to construct a shape function to be globally $C^{1}$ continuous but piecewise $C^{2}$ continuous, namely, $H^{2}$ regular, for those mesh-based numerical method. However, according to the universal approximation theorem, see Cybenko \cite{cybenko1989approximation} and Hornic \cite{hornik1991approximation}, any continuous function can be approximated arbitrarily well by a feedforward neural network, even with a single hidden layer, which offers a new possibility of analysing Kirchhoff plate bending problems. We will first give a brief introduction of deep learning.

Deep learning was first brought up as a new branch of machine learning in the realm of artificial intelligence in 2006 \cite{Vargas_2018}, which uses deep neural networks to learn features of data with high-level of abstractions \cite{lecun2015deep}. The deep neural networks adopt artificial neural network architectures with various hidden layers, which exponentially reduce the computational cost and amount of training data in some applications \cite{al2018solving}. The major two desirable traits of deep learning lie in the nonlinear processing in multiple hidden layers in supervised or unsupervised learning \cite{Vargas_2018}. Several types of deep neural networks such as convolutional neural networks (CNN) and recurrent/recursive neural networks (RNN) \cite{patterson2017deep} have been created, which further boost the application of deep learning in image processing \cite{yang2018visually}, object detection \cite{zhao2019object}, speech recognition \cite{nassif2019speech} and many other domains including genomics \cite{yue2018deep} and even finance \cite{fischer2018deep}. 

As a matter of fact, artificial neural networks (ANN) which are main tools in deep learning have been around since the 1940's \cite{mcculloch1943logical} but have not performed well until recently. They only become a major part of machine learning in the last several decades due to strides in computing techniques and explosive growth in date collection and availability, especially the arrival of backpropagation technique and advance in deep neural networks. However, based on the function approximation capabilities of feed forward neural networks, ANN were adopted to solving partial differential equations (PDEs) \cite{lagaris1998artificial,lagaris2000neural,mcfall2009artificial}, which results in a solution that can be described by a closed analytical form. Basically, ANN methods can be suitable for solving PDEs in that they are smooth enough, solutions in analytical forms can be evaluated at arbitrary points in or outside the problem domain. Yadav et al. elaborately introduced the network methods for differential equations \cite{yadav2015introduction}. In the past, when neural networks with many hidden layers were tried to solve nonlinear PDEs in order to get a better results, it usually took a long time for training, which is due to a vanishing gradient problem. However, the proposal of pretraining, which sets the initial values of connection weights and biases, with the back propagation algorithm are now proposed to solve this problem efficiently. More recently, with improved theory incorporating unsupervised pre-training, stacks of auto-encoder variants, and deep belief nets, deep learning has become a central and popular hotspot in research and applications. 

Also, some researchers studied the application of deep learning in solving PDEs. Mills et al. deployed a deep conventional neural network to solve Schr$\rm \ddot{o}$dinger equation, which directly learned the mapping between potential and energy \cite{Mills:2017aa}. E et al. applied deep learning-based numerical methods for high-dimensional parabolic PDEs and back-forward stochastic differential equations, which was proven to be efficient and accurate even for 100-dimensional nonlinear PDEs \cite{E_2017,Han:2018aa}. Also, E and Yu proposed a Deep Ritz method for solving variational problems arising from partial differential equations \cite{E_2018}. Raissi et al. however solves PDEs in a different way and has made a series of contribution to this field. They first applied the probabilistic machine learning in solving linear and nonlinear differential equations using Gaussian Processes and later introduced a data-driven Numerical Gaussian Processes to solve time-dependent and nonlinear PDEs, which circumvented the need for spatial discretization \cite{RAISSI2017683,RAISSI2018125,RaissiNGP2018}. Later, Raissi et al. \cite{Raissi:2017aa,Raissi:2017ab,RAISSI2019686} introduced a physical informed neural networks for supervised learning of nonlinear partial differential equations from Burger's equations to Navier-Stokes equations. Two distinct models were tailored for spatio-temporal datasets: continuous time and discrete time models. Raissi later employed a deep learning approach for discovering nonlinear PDEs from noisy observations in space and time with two deep neural networks, one for the representation of nonlinear-dynamic PDEs and one for a prior on the unknown solution \cite{Raissi:2018:DHP:3291125.3291150}. Raissi applied a deep neural networks in solving coupled forward-backward stochastic differential equations and their corresponding high-dimensional PDEs \cite{Raissi:2018aa}. Beck et al. \cite{Beck:2018aa,Beck_2019} studied the deep learning in solving stochastic differential equations and Kolmogorov equations, and validated the accuracy and speed proposed method, especially in high dimensions. Nabian and Meidani studied the presentation of high-dimensional random partial differential equations with a feed-forward fully-connected deep neural networks \cite{Nabian:2018aa,Nabian:2018ab}. Based on the physics informed deep neural networks, Tartakovsky et al. studied the estimation of parameters and unknown physics in PDE models \cite{Tartakovsky:2018aa}. Qin et al. applied the deep residual network and observation data to approximate unknown governing differential equations \cite{Qin:2018aa}. Sirignano and Spiliopoulos \cite{SIRIGNANO20181339} gave a theoretic motivation of using deep neural networks as PDE approximators, which converges as the number of hidden layers tend to infinity. Based on this, a deep Galerkin method was tested to solve PDEs including high-dimensional ones. Berg and Nystr$\rm \ddot{o}$m \cite{BERG201828} proposed a unified deep neural network approach to approximate solutions to PDEs and then used deep learning to discover PDEs hidden in complex data sets from measurement data \cite{BERG2019239}. In general, a deep feed-forward neural networks can well-sever as a suitable solution approximators, especially for high-dimensional PDEs with complex domains. 

Meanwhile, some researchers study the surrogate of FEM by deep learning, which mainly trains the deep neural networks from datasets obtained from FEM. From work done by Liang et al. \cite{Liang_2017,Liang_2018}, a machine learning approach was first used to investigate the relationship between geometric features of aorta and FEM-predicted ascending aortic aneurysm rupture risk and then a deep learning was used to estimate the stress distribution of the aorta, which will be beneficial to real-time patient-specific computational simulations. Lee et al. introduced the background information involved in using deep learning for structural engineering \cite{Lee_2017}. Later, Wang et al. \cite{WANG2019499} applied deep learning in calculating U* index for the high efficient load paths analysis, with training data obtained from ANSYS results. 

However, in this research, we will not confine deep learning application within FEM datasets. Rather, the deflection of Kirchhoff plate is first approximated with deep physical informed feedforward neural networks with hyperbolic tangent activation functions and trained by minimizing loss function related to governing equation of Kirchhoff bending problems and related boundary conditions. The training data for deep neural networks are obtained by randomly distributed collocation points from the physical domain of the plate. And clearly, this deep collocation method is a truly mesh-free method without the need of background grids. In this study, the method is established and applied to enrich deep learning with longstanding developments in engineering mechanics.

The paper is organised as follows: 
First a brief introduction of Kirchhoff plate bending strong form with typical boundary conditions is given. Then we introduce a basic knowledge of the deep learning technique and algorithms, which be helpful for later application. For numerical analysis, the deep collocation method with varying hidden layers and neurons are adopted for plates with various shapes, boundary and load conditions, hoping to manifest the favourable numerical features such as high accuracy and robustness of the proposed method.

\section{Kirchhoff plate bending}\label{sec2:concept}

Based on Kirchhoff plate bending theory \cite{ventsel2001thin}, the relation between lateral deflection $w\left (x,y \right )$ of the middle surface $(z=0)$ and rotations about the $x$,$y$-axis can be given by 
\begin{equation}
\theta_{x} =\frac{\partial w}{\partial x}, \qquad \theta_{y} =\frac{\partial w}{\partial y}.
\label{Rotation}
\end{equation}
\begin{figure}
\centering
\begin{tabular}{c}
\subfloat{\includegraphics[width=0.75\textwidth]{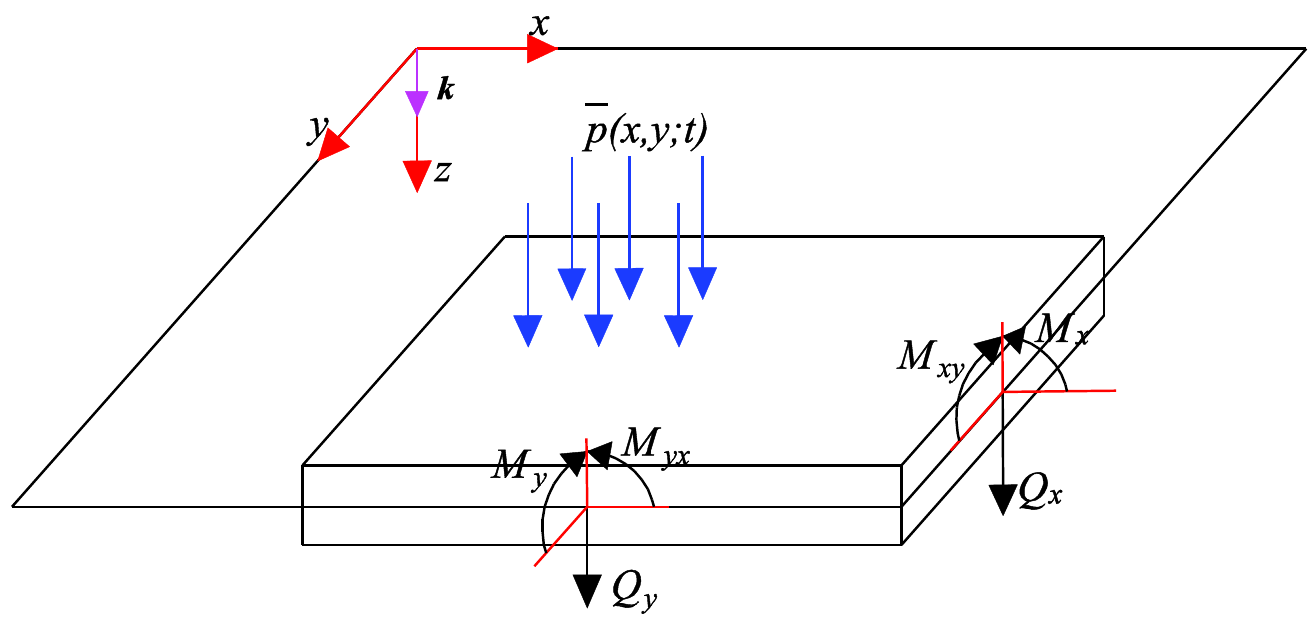}}
\end{tabular}
\caption{Kirchhoff plate in the coordinate system.}
\label{Figure1:plate}
\end{figure}
Under the coordinate system shown in Figure~\ref{Figure1:plate}, the displacement field in a thin plate can be expressed as:
\begin{equation}
\begin{aligned}
& u\left ( x,y,z \right )=-z\frac{\partial w}{\partial x}, \\
& v\left ( x,y,z \right )=-z\frac{\partial w}{\partial y}, \\
& w\left ( x,y,z \right )=w\left ( x,y \right ).
\end{aligned}
\label{Displacement field}
\end{equation}
It is obviously that the transversal deflection of the middle plane of the thin plate  can be regard as the field variables of the bending problem of thin plates. The corresponding bending and twist curvatures are the generalized strains:
\begin{equation}
	k_{x}=-\frac{\partial^2 w}{\partial x^2},\: k_{y}=-\frac{\partial^2 w}{\partial y^2},\: k_{xy}=-2\frac{\partial^2 w}{\partial x \partial y}.
\end{equation}
Therefore, the geometric equations of Kirchhoff bending can be expressed as:
\begin{equation}
\textbf{\textit{k}}=\begin{Bmatrix}
k_{xx}\\ 
k_{yy}\\ 
k_{xy}
\end{Bmatrix}=-\begin{Bmatrix}
\frac{\partial^2 w}{\partial x^2} \\[5pt]
\frac{\partial^2 w}{\partial y^2} \\[5pt]
2\frac{\partial^2 w}{\partial x \partial y}
\end{Bmatrix}=\textbf{\textit{L}}w,
\end{equation}           
with $\textbf{\textit{L}}$ being the differential operator defined as $\textbf{\textit{L}}=-\begin{pmatrix}
				\frac{\partial^2 }{\partial x^2} &\frac{\partial^2 }{\partial y^2}  & 2\frac{\partial^2 }{\partial x \partial y}
				\end{pmatrix}^{T}$. Accordingly, the bending and twisting moments, shown in Figure~\ref{Figure1:plate} can be obtained as:
\begin{equation}
\begin{aligned}
& M_{x}=-D_{0}\left ( \frac{\partial^2w }{\partial x^2} +\nu \frac{\partial^2w }{\partial y^2} \right ), \\
& M_{y}=-D_{0}\left ( \frac{\partial^2w }{\partial y^2} +\nu \frac{\partial^2w }{\partial x^2} \right ), \\
& M_{xy}=M_{yx}=-D_{0}\left ( 1-\nu \right )\frac{\partial^2w }{\partial xy}.
\end{aligned}
\label{moment}
\end{equation}				
Here $D_{0}=\frac{Eh^{3}}{12\left ( 1-\nu^{2} \right )}$	 is the bending rigidity, where $E$ and $\nu$ are the Young's modulus and Poisson ratio, and $h$ is the thickness of the thin plate.	
For isotropic thin plate, the constitutive equation	can be expressed in Matrix form
\begin{equation}
	\textbf{\textit{M=Dk}}
\end{equation}	 
with $D=D_{0}\begin{bmatrix}
1 &  \nu & 0 \\ 
  \nu& 1  &0 \\ 
 0& 0 & \left ( 1-\nu \right )/2
\end{bmatrix}$.	
The shear forces can be obtained in terms of the generalizsed stress components
\begin{equation}
Q_{x}=\frac{\partial M_{x}}{\partial x}+\frac{\partial M_{xy}}{\partial y},\, \, Q_{y}=\frac{\partial M_{xy}}{\partial x}+\frac{\partial M_{y}}{\partial y}
\label{shearforce}
\end{equation}	

The differential equation for the deflections for thin plate based on Kirchhoff's assumptions can be expressed by transversal deflection as
\begin{equation}
\bigtriangledown^{2}\left ( \bigtriangledown^{2}w \right )=\bigtriangledown^{4}w=\frac{p}{D}
\label{governing}
\end{equation}
where $\bigtriangledown^{4}\left (  \right )=\frac{\partial^4 }{\partial x^4}+2\frac{\partial^4 }{\partial x^2\partial y^2}+\frac{\partial^4 }{\partial y^4}$ is commonly called biharmonic operator.	

Consequently, the Kirchhoff plate bending problems can be boiled down to a fourth order PDE problem, which pose difficulty for tradition mesh-based method in constructing a shape function to be 	$H^{2}$ regular. 
Moreover, the boundary conditions of Kirchhoff plate taken into consideration in this paper can be generally classified into three parts, namely,
\begin{equation}
\partial\Omega =\Gamma_{1}+\Gamma_{2}+\Gamma_{3}.
\label{boundary}
\end{equation}

For clamped edge boundary, $\Gamma_{1}:  w=\tilde{w}, \  \frac{\partial w}{\partial n} = \tilde{ \theta}_{n} $, $w=\tilde{w}, \  \tilde{ \theta}_{n} $ are functions of arc length along this boundary.
 
For simply supported edge boundary, $\Gamma_{2}:  w=\tilde{w}, \  M_{n} =\tilde{ M}_{n}$, $\tilde{ M}_{n} $ is also a function of arc length along this boundary. 

For free boundary conditions, $\Gamma_{3}:  M_{n} =\tilde{ M}_{n} , \  \frac{\partial M_{ns} }{\partial s}+Q_{n}=\tilde{q}$, where $\tilde{q}$ is the load exerted along this boundary.

It should be noted that $\textit{\textbf{n}},\textit{\textbf{s}}$ here refer to the normal and tangent directions along the boundaries.

\section{Deep Collocation Method for solving Kirchhoff plate bending}\label{sec3:methodology}
In this section, we will begin with introducing some preliminaries on deep learning, including the feed forward neural network architectures, some useful algorithms involved in deep learning. Then based on those basis, the formulation of deep collocation method is elucidated. 

\subsection{Feed forward neural network}
The basic architecture of a fully connected feedforward neural network is shown in Figure~\ref{Figure2:network}, which comprises of multiple layers: input layer, one or more hidden layers and output layer. Each layer consists of one or more nodes called neurons, shown in the Figure~\ref{Figure2:network} by small coloured circles, which is the basic unit of computation. For an interconnected structure, every two neurons in neighbouring layers have a connection, which is represented by a connection weight. Depicted in Figure~\ref{Figure2:network}, the weight between neuron $k$ in hidden layer $l-1$ and neuron $j$ in hidden layer $l$ is denoted by $w_{jk}^{l}$. No connection exists among neurons in the same layer as well as in the non-neighbouring layers. Input data, defined from $x_{1}$ to $x_{N}$, flow through this neural network via connections between neurons, starting from input layer, through hidden layer $l-1$, $l$, to output layer, which eventually output data from $y_{1}$ to $y_{M}$. The feedforward neural work defines a mapping $FNN: \mathbb{R}^N \to \mathbb{R}^M$.
\tikzset{
  every neuron/.style={
    circle,
    draw,
    minimum size=1.2cm,
  },
  neuron missing/.style={
    draw=none,
    scale=2,
    text height=0.25cm,
    fill=none,
    execute at begin node=\color{black}$\vdots$
  },
}

\begin{figure}
\centering
\begin{tikzpicture}[x=1.75cm, y=1.2cm, >=stealth]

\tikzstyle{input neuron}=[every neuron, fill=green!50]
\tikzstyle{output neuron}=[every neuron, fill=red!50]
\tikzstyle{hidden neuron}=[every neuron, fill=blue!50]

\node [] at (3,3) {$Forward\:propagation\:of\:activation\:values$};
\draw[->,semithick] (0,2.7) -- (6,2.7);

\foreach \m/\l [count=\y] in {1,missing,2}
	\node [input neuron/.try, neuron \m/.try] (input-\m) at (0,1.5-\y) {};

\foreach \m [count=\y] in {1,missing}
	\node [hidden neuron/.try, neuron \m/.try ] (hidden1-\m) at (2,2.5-\y) {};
\foreach \m [count=\y] in {2}
	\node [hidden neuron/.try, neuron \m/.try ] (hidden1-\m) at (2,2.5-\y-2) {$b^{l-1}_k$};
\foreach \m [count=\y] in {missing,3}
	\node [hidden neuron/.try, neuron \m/.try ] (hidden1-\m) at (2,2.5-\y-3) {};

\foreach \m [count=\y] in {1,missing}
	\node [hidden neuron/.try, neuron \m/.try ] (hidden2-\m) at (4,2.5-\y) {};
\foreach \m [count=\y] in {2}
	\node [hidden neuron/.try, neuron \m/.try ] (hidden2-\m) at (4,2.5-\y-2) {$b^l_j$};
\foreach \m [count=\y] in {missing,3}
	\node [hidden neuron/.try, neuron \m/.try ] (hidden2-\m) at (4,2.5-\y-3) {};

\foreach \m [count=\y] in {1,missing,2}
  \node [output neuron/.try, neuron \m/.try ] (output-\m) at (6,1.5-\y) {};

\foreach \l [count=\i] in {1,N}
    \node [] at (input-\i) {$x_{\l}$};

\node [] at (output-1) {$b_1$};
\node [] at (output-2) {$b_M$};

\foreach \l [count=\i] in {1,M}
  \draw [->] (output-\i) -- ++(1,0)
    node [above,midway] {$y_\l$};

\foreach \i in {1,...,2}
  \foreach \j in {1,...,3}
    \draw [dashed,->] (input-\i) -- (hidden1-\j);

\foreach \i in {1}
  \foreach \j in {1,...,3}
    \draw [->] (hidden1-\i) -- (hidden2-\j);

\foreach \i in {2}
  \foreach \j in {1}
    \draw [->] (hidden1-\i) -- (hidden2-\j);

\foreach \i in {2}
  \foreach \j in {2}
    \draw [->] (hidden1-\i) -- (hidden2-\j)
    	node [above,midway] {$w^l_{jk}$};

\foreach \i in {2}
  \foreach \j in {3}
    \draw [->] (hidden1-\i) -- (hidden2-\j);

\foreach \i in {3}
  \foreach \j in {1,...,3}
    \draw [->] (hidden1-\i) -- (hidden2-\j);

\foreach \i in {1,...,3}
  \foreach \j in {1,...,2}
    \draw [dashed,->] (hidden2-\i) -- (output-\j);

\foreach \l [count=\x from 0] in {Input Layer, Hidden Layer $l-1$, Hidden Layer $l$, Output Layer}
  \node [align=center, above] at (\x*2,2) {\l};

\node [] at (3,-3.3) {$Back\:propagation\:of\:errors$};
\draw[<-,semithick] (0,-3.6) -- (6,-3.6);

\end{tikzpicture}
\caption{Architecture of a fully connected feedforward back-propagation neural network.}
\label{Figure2:network}
\end{figure}
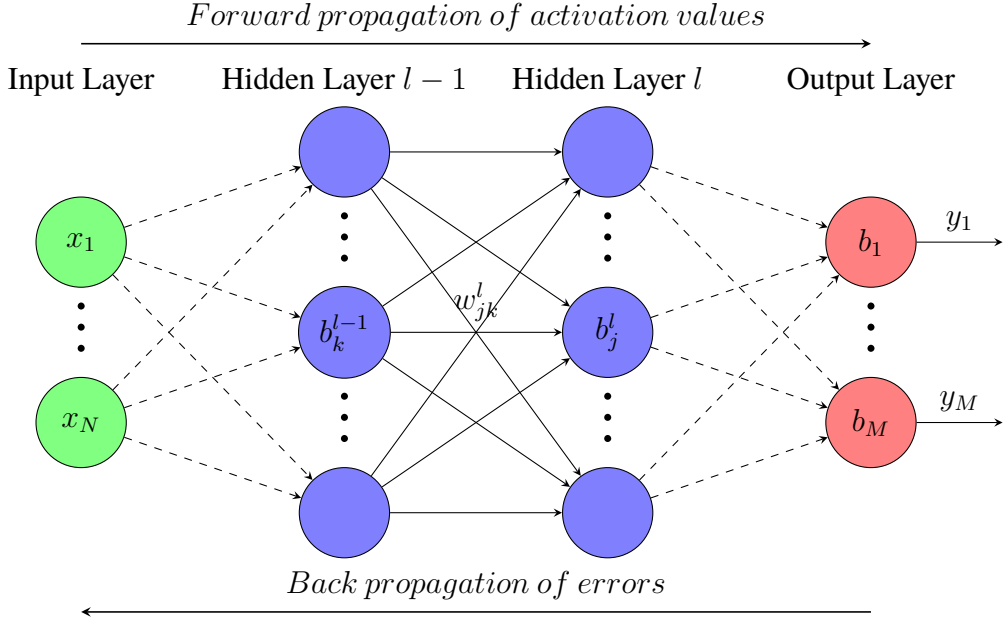

However, it should be noted that the number of neurons on each hidden layers and number of hidden layers can be any number and are invariably determined through a trial and error procedure. It has also been concluded that any continuous function can be approximated with any desired precision by a feed forward with even a single hidden layer \cite{FUNAHASHI1989183,HORNIK1989359}.

On each neuron in the feed forward neural network, a bias is supplied including neurons in the output layer except the neurons in the input layer, which is defined by $b^l_j$ for bias of neuron $j$ in layer $l$. Besides, the activation function is defined for output of each neuron in order to introduce a non-linearity into the neural network and make the back-propagation possible where gradients are supplied along with an error to update weights and biases. The activation function in layer $l$ will be denoted by $\sigma$ here. There are many activation functions can be used such as sigmoids function, hyperbolic tangent function $ \left( Tanh \right)$, Rectified linear units $ \left( Relu \right)$, and so on. Some suggestions upon the choice of activation function can be referred in \cite{Hayou:2018aa}. Hence, for the value on each neuron in the hidden layers and output layer adds the weighted sum of values of output values from the previous layer with corresponding connection weights to basis on the neuron. A intermediate quantity for neuron $j$ on hidden layer $l$ is defined as
\begin{equation}
a^l_j = \sum_k w^l_{jk}y^{l-1}_k + b^l_j,
\label{activation}
\end{equation}
and its output is given by the activation of the above weighted input
\begin{equation}
y^{l}_j =\sigma \left( a^l_j \right)=\sigma \left(\sum_k w^l_{jk}y^{l-1}_k + b^l_j \right),
\label{wi} 
\end{equation}
where $y^{l-1}_k$ is the output from previous layer. 

So, basically, when Equation \ref{wi} is applied to compute $y^{l}_j$, the intermediate quantity $a^l_j$ was calculated along the way. This quantity turns out to be useful and named here as weighted input to neuron $j$ on hidden layer $l$. Equation \ref{activation} can be written in a compact matrix form, which calculate weighted inputs for all neurons on certain layer efficiently, obtaining:
\begin{equation}
\textit{\textbf{a}}= \textit{\textbf{W}}^{l}\textit{\textbf{y}}^{l-1}+\textit{\textbf{b}}^{l},
\label{a} 
\end{equation}
and accordingly, from Equation \ref{a}, $\textit{\textbf{y}}^{l}=\sigma\left( \textit{\textbf{a}} \right)$, where activation function is applied elementwise. A feedforward network thus defines a function $\textit{f}\left(\textit{\textbf{x}};\theta\right)$ depending on input data $\textit{\textbf{x}}$ and parametrised by $\theta$ consisting of weights and biases in each layer. The defined function provides an efficient way to approximate unknown field variables.

\subsection{Backpropagation}
Backpropagation $\left( backward\;propagation \right)$ is an important and computationally efficient mathematical tool to compute gradients in deep learning \cite{nielsenneural}. Essentially, backpropagation is based on recursively applying the chain rule and decides which computations can be run in parallel from computational graphs. In our problem, the governing equation is the fourth order partial derivatives of field variable $w\left(\textit{\textbf{x}}\right)$ approximated by the deep neural networks $\textit{f}\left(\textit{\textbf{x}};\theta\right)$, so this makes backpropagation a critical role. 
For the approximation defined by $\textit{f}\left(\textit{\textbf{x}};\theta\right)$, in order to find the weights and biases, a loss function $\textit{\textbf{L}}\left(\textit{f},w\right)$ is defined to be minimised \cite{Janocha_2017}. The backpropagation algorithm for computing the gradient of loss function $\textit{\textbf{L}}\left(\textit{\textbf{f}},w\right)$ can be defined as follows \cite{nielsenneural}:
\begin{itemize}
\item \textbf{Input}: Input dataset $x^{1},...,x^{n}$, prepare activation $y^1$ for input layer;
\item \textbf{Feedforward}: For each layer $xl=2,3,...,L$, compute $a^l = \sum_k W^ly^{l-1} + b^l$, and $\sigma \left( a^l \right)$;
\item \textbf{Output error}: Compute the error $\delta^L= \nabla_{y^L} \textit{\textbf{L}} \odot \sigma'_L(a^L)$
\item \textbf{Backpropagation error}: For each $l=L-1,L-2,...,2$, compute $\delta^l = \left((W^{l+1})^T \delta^{l+1}\right) \odot \sigma'_l(a^l)$;
\item \textbf{Output}: The gradient of the loss function is given by $\frac{\partial \textit{\textbf{L}}}{\partial w^l_{jk}} = y^{l-1}_k \delta^l_j$ and $\frac{\partial \textit{\textbf{L}}}{\partial b^l_j} = \delta^l_j$.
\end{itemize}
 Here, $\odot$ denotes the Hadamard product.
 
Now, there are a lists of deep learning frameworks for us to choose to setup a training. The main two approaches, Pytorch and Tensorflow, however computing derivatives in the computational graphs distinctly. The former inputs a numerical value and then compute the derivatives at this node, while the latter computers the derivatives of a symbolic variable, then store the derivative operations into new nodes added to the graph for later use. Obviously, the latter is more advantageous in computing higher-order derivatives, which can be computed from its extended graph by running backpropagation repeatedly. In this paper, since the fourth-order derivatives of field variables is needed to be computed, the Tensorflow framework is thus adopted for calculation \cite{Al-Aradi:2018aa}.

\subsection{Formulation of deep collocation method}
The formulation of a deep collocation in solving Kirchhoff plate bending problems is introduction in this section. Collocation method is a widely used method seeking numerical solutions for ordinary, partial differential and integral equations \cite{atluri2005methods}. It is a popular method for trajectory optimization in control theory. A set of randomly distributed points (also known as collocation points) is often deployed to represent a desired trajectory that minimizes the loss function while satisfying a set of constraints. The collocation methods tend to be relatively insensitive to instability of system (such as blowing/vanishing gradients with neural networks), then it can be a viable way to train the deep neural networks in this paper \cite{agrawalcollocation}.

Recalled form Section 2, Equation \ref{governing},\ref{boundary}, the solving of Kirchhoff plate bending problems can be boiled down to the solving of a fourth order biharmonic equations with the type of boundary constraints. Thus we first discretize the physical domain with collocation points denoted by $\textit{\textbf{x}}\,_\Omega=(x_1,...,x_{N_\Omega})^T$. Another set of collocation points are deployed to discretize boundary conditions denoted by $\textit{\textbf{x}}\,_\Gamma(x_1,...,x_{N_\Gamma})^T$.
Then the transversal deflection $w$ is approximated with the aforementioned deep feedforward neural network $\textit{w}^h \left(\textit{\textbf{x}};\theta\right)$. A loss function can thus be constructed to find the approximate solution by considering the minimizing of governing equation with boundary conditions approximated by $\textit{w}^h \left(\textit{\textbf{x}};\theta\right)$. The mean squared error loss form is adopted here.

Substituting $\textit{w}^h \left(\textit{\textbf{x}}\,_\Omega;\theta\right)$ into Equation \ref{governing}, we can get:
\begin{equation}
G\left(\textit{\textbf{x}}\,_\Omega;\theta\right)=\bigtriangledown^{4}w^h\left(\textit{\textbf{x}}\,_\Omega;\theta\right)-\frac{p}{D},
\end{equation}
which results in a physical informed deep neural network $G\left(\textit{\textbf{x}}\,_\Omega;\theta\right)$. 

For boundary conditions illustrated in Section 2, considering all three boundaries, they can also be expressed by the neural network approximation $\textit{w}^h \left(\textit{\textbf{x}}\,_\Gamma;\theta\right)$ as:

\noindent On $\Gamma_{1}$, we have
\begin{equation}
  \textit{w}^h \left(\textit{\textbf{x}}\,_{\Gamma_1};\theta\right)=\tilde{w}, \  \frac{\partial \textit{w}^h \left(\textit{\textbf{x}}\,_{\Gamma_1};\theta\right)}{\partial n} = \tilde{ \theta}_{n}.
\end{equation} 

\noindent On $\Gamma_{2}$,
\begin{equation}
  \textit{w}^h \left(\textit{\textbf{x}}\,_{\Gamma_2};\theta\right)=\tilde{w}, \ \tilde{ M}_{n}\left(\textit{\textbf{x}}\,_{\Gamma_2};\theta\right)=\tilde{ M}_{n}, 
\label{bd2}
\end{equation} 
where $\tilde{ M}_{n}\left(\textit{\textbf{x}}\,_{\Gamma_2};\theta\right)$ can be obtained from Equation \ref{moment} by combing $\textit{w}^h \left(\textit{\textbf{x}}\,_{\Gamma_2};\theta\right)$.

\noindent On $\Gamma_{3}$,
\begin{equation}
  M_{n}\left(\textit{\textbf{x}}\,_{\Gamma_3};\theta\right) =\tilde{ M}_{n} , \  \frac{\partial M_{ns}\left(\textit{\textbf{x}}\,_{\Gamma_3};\theta\right) }{\partial s}+Q_{n}\left(\textit{\textbf{x}}\,_{\Gamma_3};\theta\right)=\tilde{q}, 
\end{equation}  
where $M_{ns}\left(\textit{\textbf{x}}\,_{\Gamma_3};\theta\right)$ can be obtained from Equation \ref{moment} and $Q_{n}\left(\textit{\textbf{x}}\,_{\Gamma_3};\theta\right)$ can be obtained from Equation \ref{shearforce} by combing $\textit{w}^h \left(\textit{\textbf{x}}\,_{\Gamma_3};\theta\right)$.

It should be noted that $\textit{\textbf{n}},\textit{\textbf{s}}$ here refer to the normal and tangent directions along the boundaries. As induced physical informed neural network $G\left(\textit{\textbf{x}};\theta\right)$, $M_{n}\left(\textit{\textbf{x}};\theta\right)$, $M_{ns}\left(\textit{\textbf{x}};\theta\right)$, $Q_{n}\left(\textit{\textbf{x}};\theta\right)$  share the same parameters as $\textit{w}^h \left(\textit{\textbf{x}};\theta\right)$. Considering the generated collocation points in domain and on boundaries, they can all be learned by minimizing the mean square error loss function:
\begin{equation}
L\left(\theta\right)=MSE=MSE_{G}+MSE_{\Gamma_{1}}+MSE_{\Gamma_{2}}+MSE_{\Gamma_{3}},
\end{equation}
with
\begin{equation}
\begin{aligned}
&MSE_{G}=\frac{1}{N_d}\sum_{i=1}^{N_d}\begin{Vmatrix}
G\left(\textit{\textbf{x}}\,_\Omega;\theta\right)
\end{Vmatrix}^2=\frac{1}{N_\Omega}\sum_{i=1}^{N_\Omega}\begin{Vmatrix}
\bigtriangledown^{4}w^h\left(\textit{\textbf{x}}\,_\Omega;\theta\right)-\frac{p}{D}
\end{Vmatrix}^2,\\
&MSE_{\Gamma_{1}}=\frac{1}{N_{\Gamma_1}}\sum_{i=1}^{N_{\Gamma_1}}\begin{Vmatrix}
 \textit{w}^h \left(\textit{\textbf{x}}\,_{\Gamma_1};\theta\right)-\tilde{w}
\end{Vmatrix}^2+\frac{1}{N_{\Gamma_1}}\sum_{i=1}^{N_{\Gamma_1}}\begin{Vmatrix}
\frac{\partial \textit{w}^h \left(\textit{\textbf{x}}\,_{\Gamma_1};\theta\right)}{\partial n} - \tilde{ \theta}_{n}
\end{Vmatrix}^2,\\
&MSE_{\Gamma_{2}}=\frac{1}{N_{\Gamma_2}}\sum_{i=1}^{N_{\Gamma_2}}\begin{Vmatrix}
 \textit{w}^h \left(\textit{\textbf{x}}\,_{\Gamma_2};\theta\right)-\tilde{w}
\end{Vmatrix}^2+\frac{1}{N_{\Gamma_2}}\sum_{i=1}^{N_{\Gamma_2}}\begin{Vmatrix}
\tilde{ M}_{n}\left(\textit{\textbf{x}}\,_{\Gamma_2};\theta\right)-\tilde{ M}_{n}
\end{Vmatrix}^2,\\
&MSE_{\Gamma_{3}}=\frac{1}{N_{\Gamma_3}}\sum_{i=1}^{N_{\Gamma_3}}\begin{Vmatrix}
 \tilde{ M}_{n}\left(\textit{\textbf{x}}\,_{\Gamma_3};\theta\right)-\tilde{ M}_{n}
\end{Vmatrix}^2+\frac{1}{N_{\Gamma_3}}\sum_{i=1}^{N_{\Gamma_3}}\begin{Vmatrix}
\frac{\partial M_{ns}\left(\textit{\textbf{x}}\,_{\Gamma_3};\theta\right) }{\partial s}+Q_{n}\left(\textit{\textbf{x}}\,_{\Gamma_3};\theta\right)-\tilde{q}
\end{Vmatrix}^2,\\
\end{aligned}
\end{equation}
where $x\,_\Omega \in {R^N} $, $\theta \in {R^K}$ is the neural network parameters. If $L\left(\theta\right)=0$, $\textit{w}^h \left(\textit{\textbf{x}};\theta\right)$ is then a solution to transversal deflection.
Our goal becomes to find the a set of parameters $\theta$ that the  approximated deflection $\textit{w}^h \left(\textit{\textbf{x}};\theta\right)$ minimize the loss $L\left(\theta\right)$. And if $L\left(\theta\right)$ is a very small value, then the approximation $\textit{w}^h \left(\textit{\textbf{x}};\theta\right)$ is very closely satisfying governing equations and boundary conditions, namely
\begin{equation}
\textit{w}^h = \argmin_{\theta \in R^K} L\left(\theta\right)
\end{equation}

Then, the solving of thin plate bending problems by deep collocation method can be reduced to an optimization problem. In deep learning Tensorflow/Pytorch framework, there are a variety available optimizers. One of the most widely used optimization method can be gradient descent based method is the Adam optimization algorithm \cite{Kingma2015AdamAM}, which is also adopted in the numerical study in this paper. Take a descent step at collocation point $\textit{\textbf{x}}_{i}$ with Adam-based learning rates $\alpha_i$, 
\begin{equation}
\theta_{i+1} = \theta_{i} + \alpha_i \bigtriangledown_{\theta } L \left ( \textit{\textbf{x}}_i;\theta_i \right )
\label{Adma}
\end{equation}
And then the process in Equation \ref{Adma} is repeated until convergence criterion is satisfied.

\section{Numerical examples}\label{sec4:example}
In this section, several numerical examples on plate bending problems with various shapes and boundary conditions is studied. And for implementation, a combined optimizer suggested by Berg et al. in \cite{BERG201828} is adopted using L-BFGS optimizer \cite{Liu_1989} first and in linear search where BFGS may fail, a  Adam optimizer is then applied with a very small learning rate. For all numerical examples, predicted maximum transverse with increasing layers are studied in order to show a convergence of deep collocation method in solving the plate bending problem.

\subsection{Simply-supported square plate}
\label{section 1:Simply-supported thin plate} 
\begin{figure}[H]
\centering
\begin{tabular}{c}
\subfloat{\includegraphics{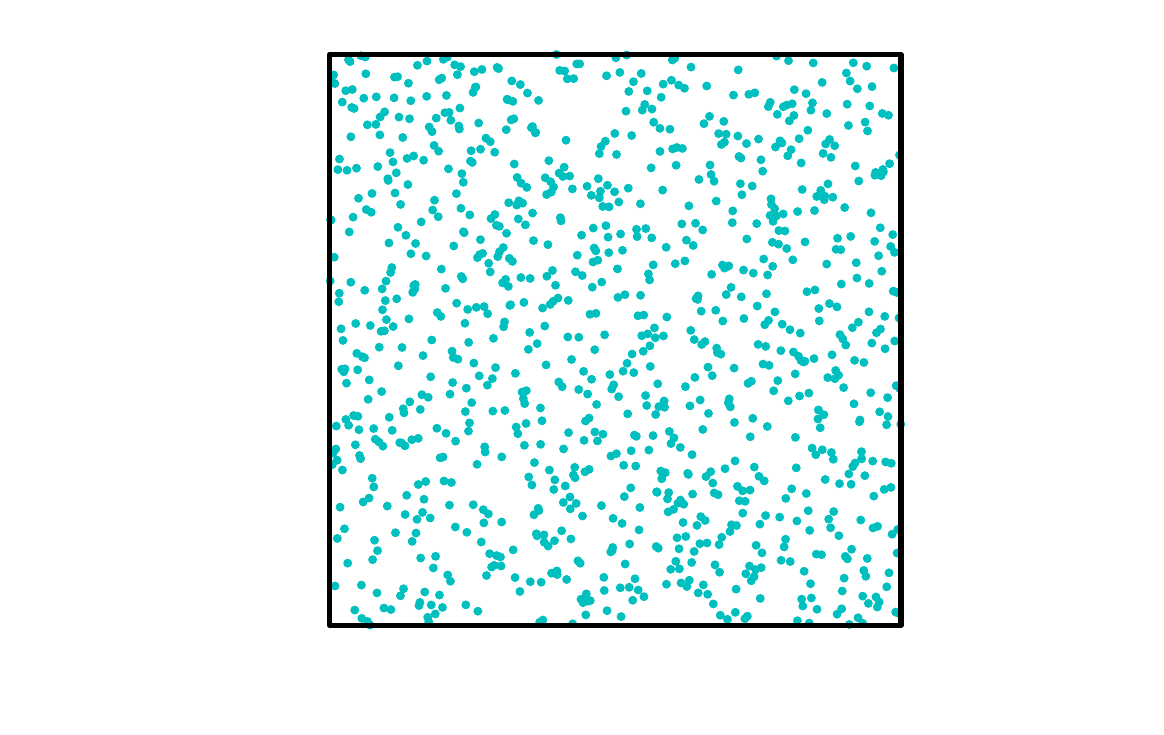}}
\end{tabular}
\vspace{-1.25cm}
\caption{Collocation points discretize the square domain.}
\label{Figure3:Scatterpoint}
\end{figure}

A simply-supported square plate under a sinusoidal distribution of transverse loading is studied. The distributed load is given by
\begin{equation}
\begin{array}{l}
p=\frac{p_{0}}{D}\textrm{sin}\left (\frac{\pi x}{a}  \right )\textrm{sin}\left (\frac{\pi y}{b}  \right ).
\end{array}
\end{equation} 
Here, $a$,$b$ the length of the plate. $D$ denotes the flexural stiffness of the plate
   and depends on the plate thickness and material properties.The exact solution for this problem is given by 
   
\begin{equation}
\begin{array}{l}
w=\frac{p_{0}}{\pi^{4} D\left (\frac{1}{a^{2}}+\frac{1}{b^{2}}  \right )^{2}}\textrm{sin}\left (\frac{\pi x}{a}  \right )\textrm{sin}\left (\frac{\pi y}{b}  \right ).
\end{array}
\end{equation}   
Here, $w$ represents the transverse plate deflection. 
For this numerical example, we first generate 1000 randomly distributed collocation points in the physical domain depicted in Figure~\ref{Figure3:Scatterpoint}. And we thoroughly studied the influence of deep neural network with a varying number of hidden layer and neurons on the maximum deflection at the centre of the plate, which is then shown in Table 1. The numerical results are compared with the exact solution. It is clear that the results predicted by more hidden layers are more desirable, especially for neural networks with three hidden layers. To better reflect the deflection vector in the whole physical domain, the contour plot, contour error plot of deflection for increasing hidden layers with 50 neurons are shown in Figure~\ref{Figure5:ssneuron50point1}, Figure~\ref{Figure6:ssneuron50point2}, Figure~\ref{Figure7:ssneuron50point3}

\begin{figure}[H]
\centering
\caption*{Table 1: Maximum deflection predicted by deep collocation method.}
\vspace{-0.55cm}
\begin{tabular}{c}
\subfloat{\includegraphics[width=0.85\textwidth]{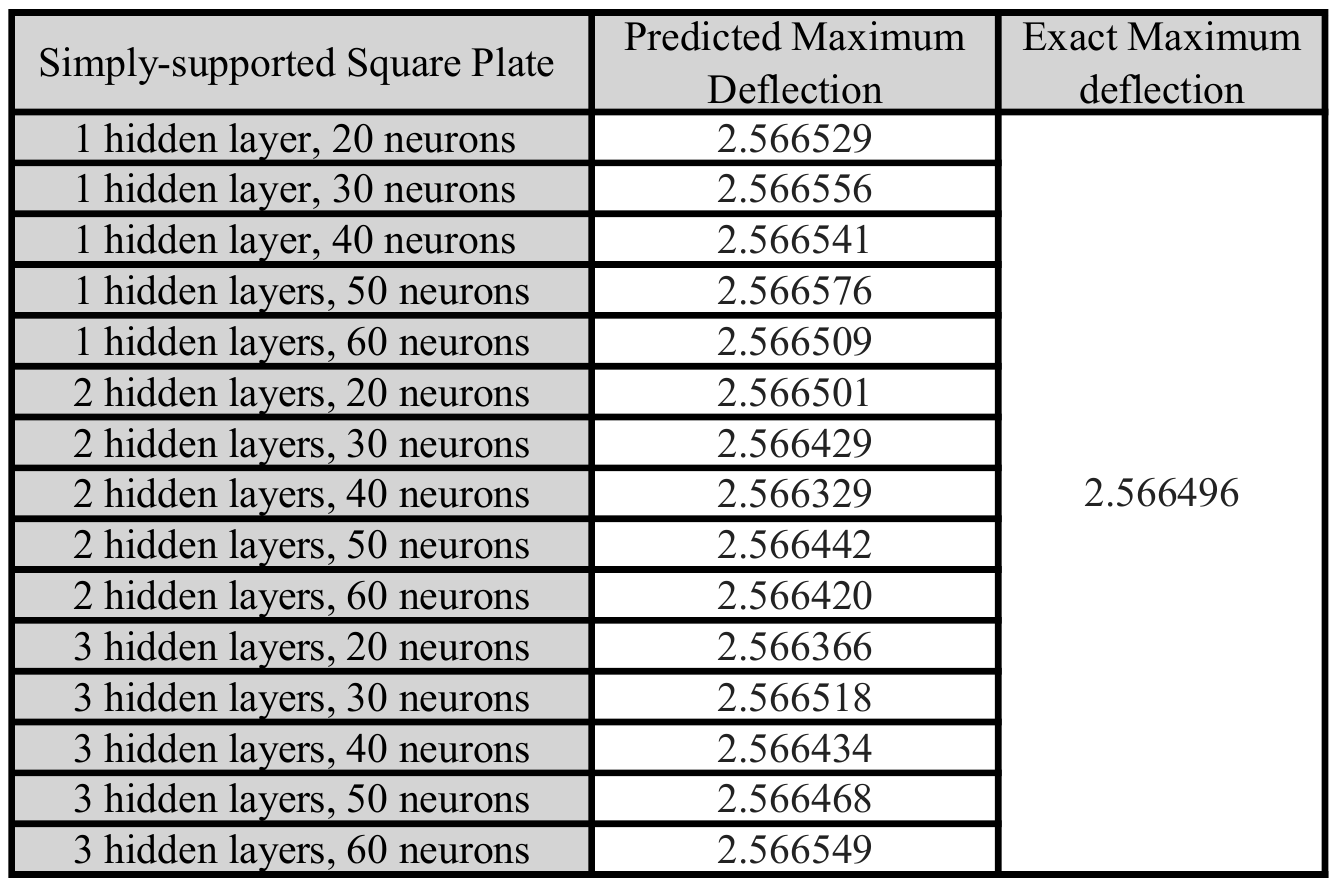}}
\end{tabular}
\end{figure}

\begin{figure}[H]
\centering
\begin{tabular}{c}
\subfloat{\includegraphics{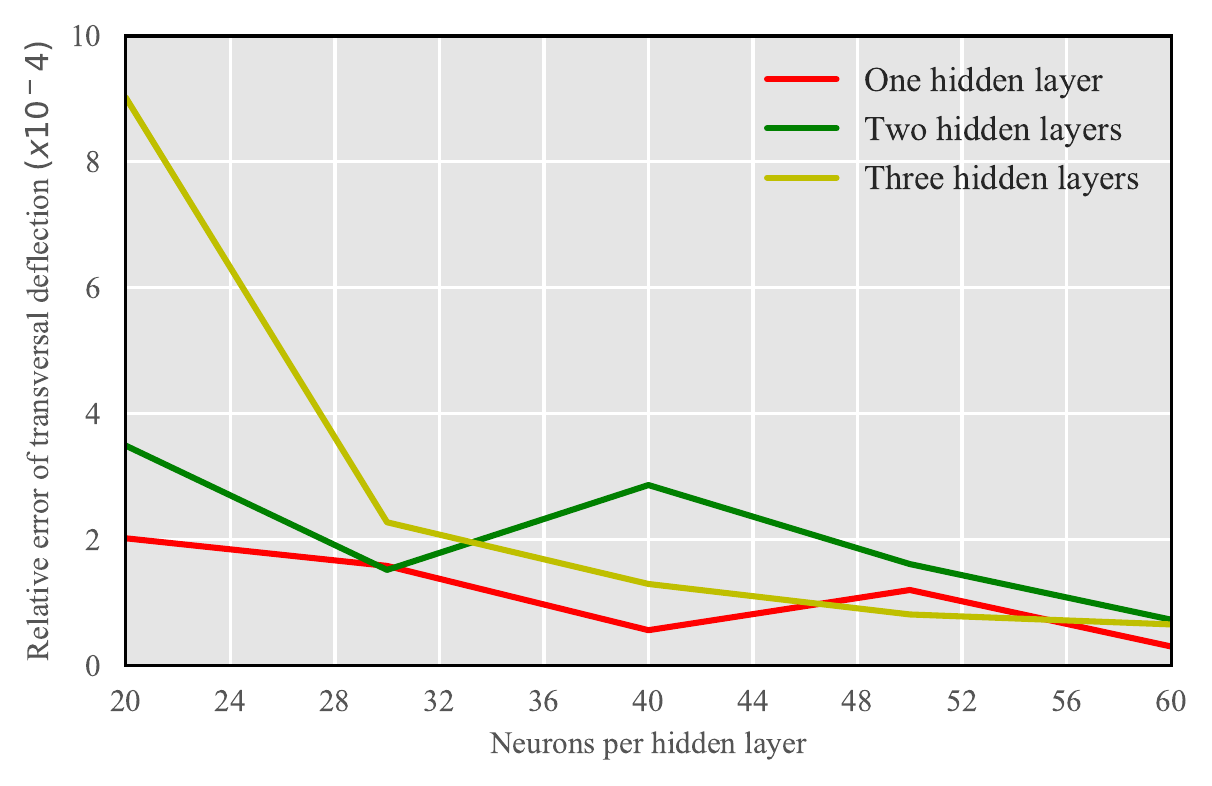}}
\end{tabular}
\vspace{-0.55cm}
\caption{The relative error of deflection with varying hidden layers and neurons.}
\label{Figure4:Relativeerrorss}
\end{figure}

In Table 1, we employed a varying number of hidden layers from 1 to 4 and in each layer and the number of neurons varies from 20 to 60. We calculated the corresponding maximum transversal deflection at the centre of the square plate. From the $L_2$ relative error of deflection vector at all predicted points is shown in Figure~\ref{Figure4:Relativeerrorss} for each case. And it is very clear for even the neural network with only one single hidden layer with 20 neurons, the results is already very accurate and favourable. For most cases, with increasing neurons and hidden layers, the results converge to the exact solution and the results are very accurate even with a few neurons and a single hidden layer. In Figure~\ref{Figure4:Relativeerrorss}, all three hidden layer types get very accurate results. Though the single layer with 20 neurons is the most accurate in all three types with 20 neurons, the magnitude of all is $1\times10^{-4}$ and the other two results are also very accurate. And as the number of hidden layer and neurons increases, the relative error curves become flat and obtain results around exact solutions.

From Figure~\ref{Figure5:ssneuron50point1}, Figure~\ref{Figure6:ssneuron50point2}, Figure~\ref{Figure7:ssneuron50point3}, we can observe that the deflection is accurately predicted by the deep collocation method, which agree well with the exact solutions.
And as the hidden layer number increases, the numerical results converge to the exact solutions in the whole square plate. The predicted plate deformation agrees well with the exact deformation. All these lend some credence to the suitable application of this deep learning based method. The advantageous of neural networks with hidden layers is not conspicuously reflected in this numerical example, as the next numerical example shows more clearly.
\begin{figure}[H]
\centering
\begin{tabular}{cc}
\subfloat[Predicted deflection contour]{\includegraphics[width=7cm,height=6cm]{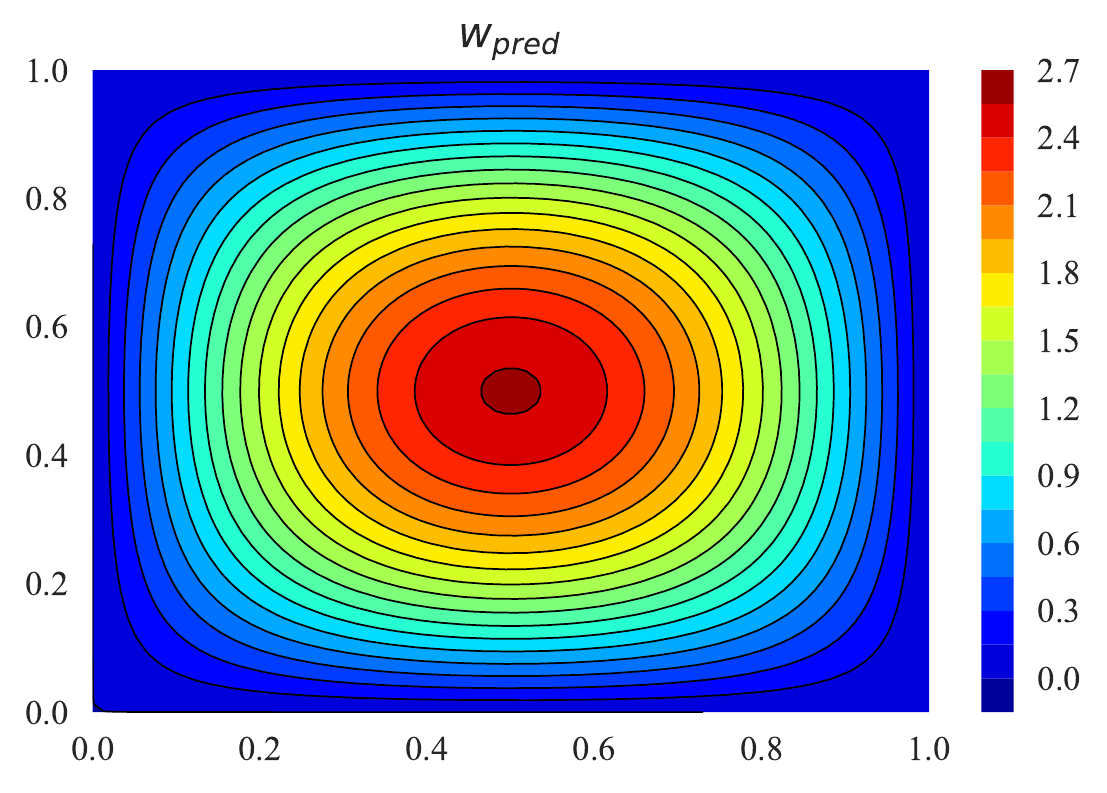}} & 
\hspace{-0.75cm}
\subfloat[Deflection error contour]{\includegraphics[width=7cm,height=6cm]{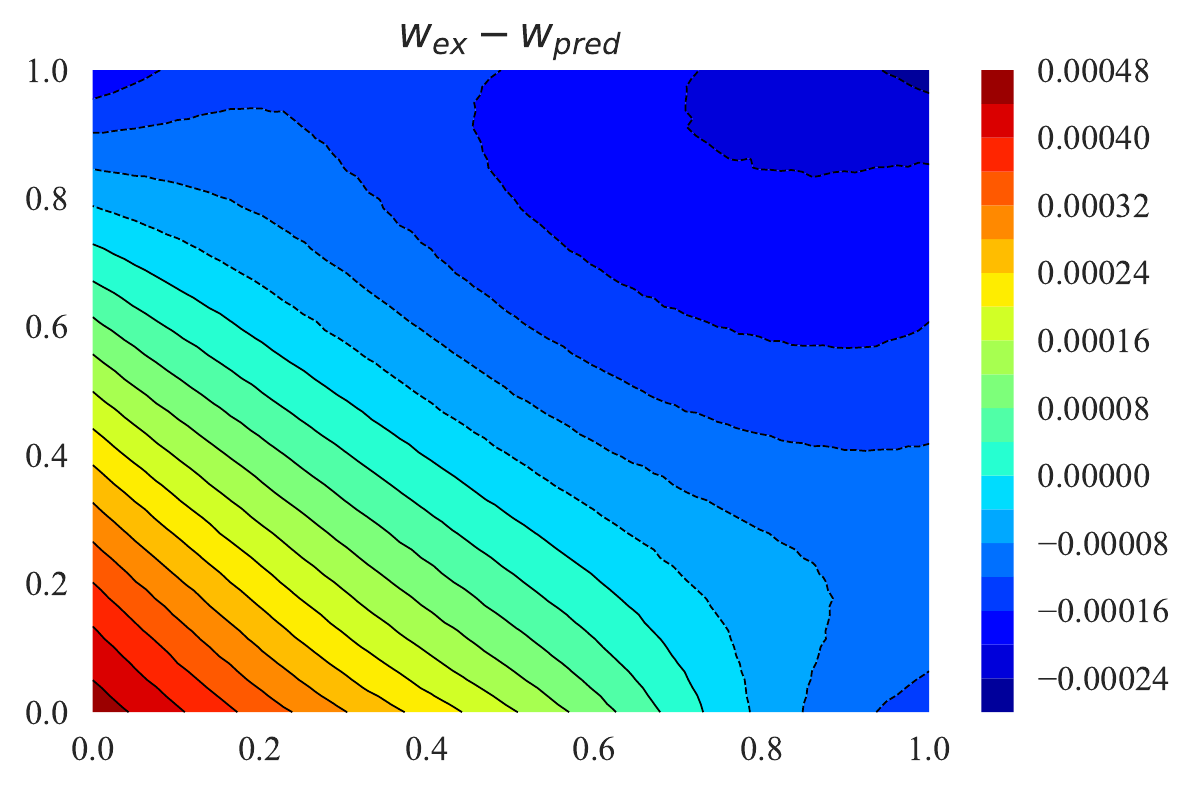}} \\ 
\vspace{-3cm}
\subfloat[Predicted deflection]{\includegraphics[width=7cm,height=6cm]{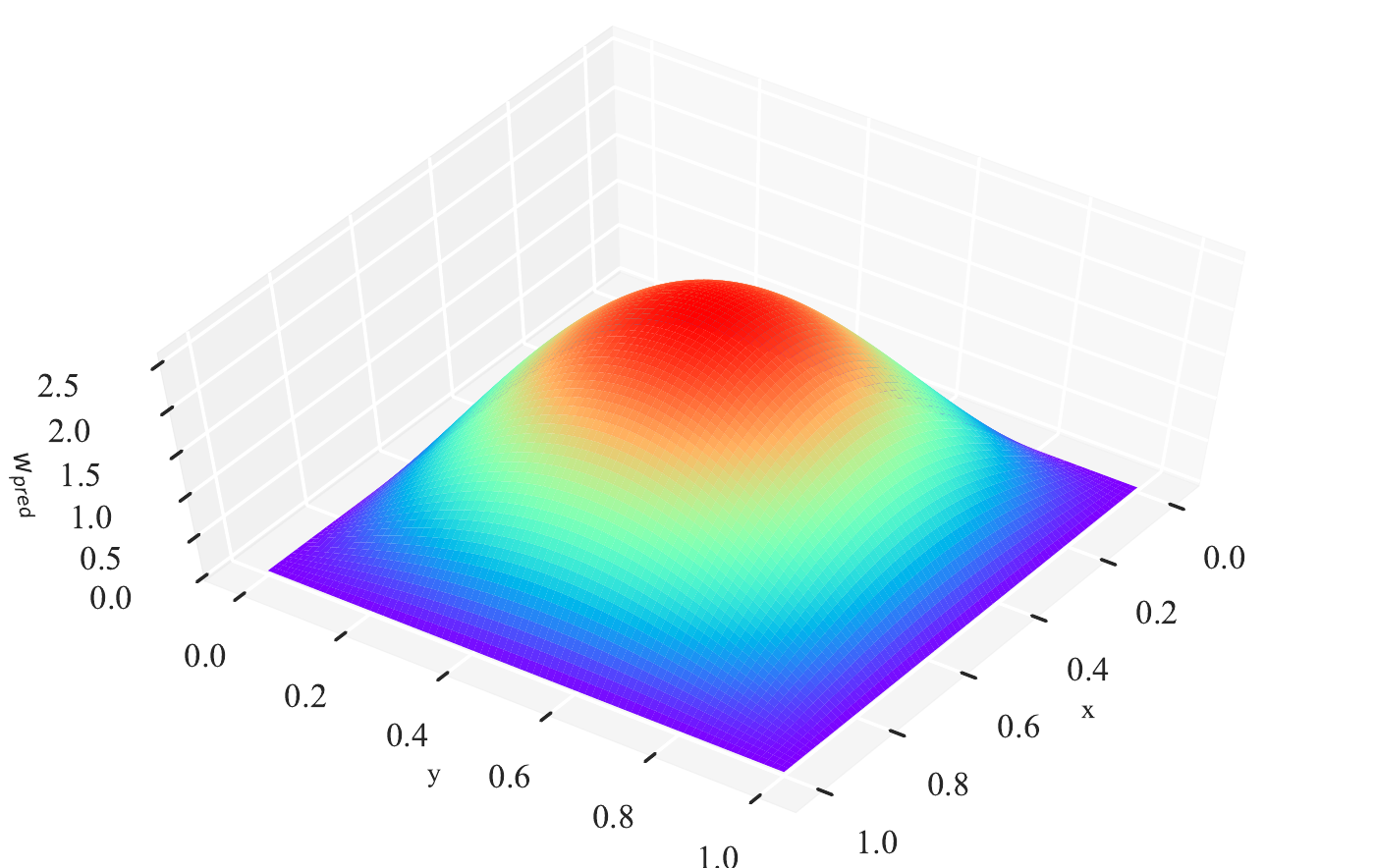}}&
\hspace{-0.75cm}
\subfloat[Exact deflection]{\includegraphics[width=7cm,height=6cm]{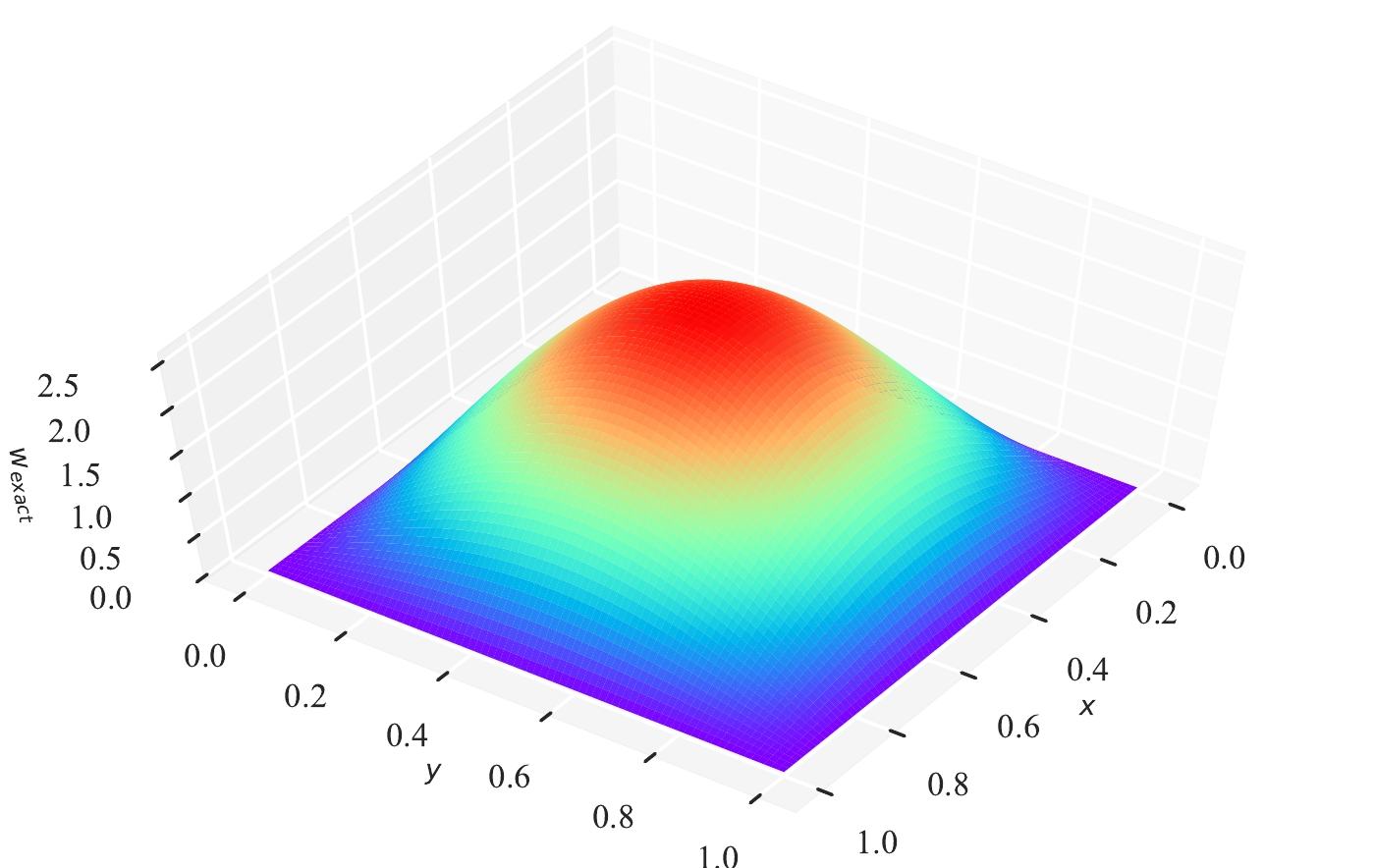}} \\
\end{tabular}
\vspace{3cm}
\caption{$\left(a\right)$ Predicted deflection contour $\left(b\right)$ Deflection error contour $\left(c \right)$ Predicted deflection $\left(d \right)$ Exact deflection of the simply-supported square plate with 1 hidden layers and 50 neurons.}
\label{Figure5:ssneuron50point1}
\end{figure}

\begin{figure}[H]
\centering
\begin{tabular}{cc}
\subfloat[Predicted deflection contour]{\includegraphics[width=7cm,height=6cm]{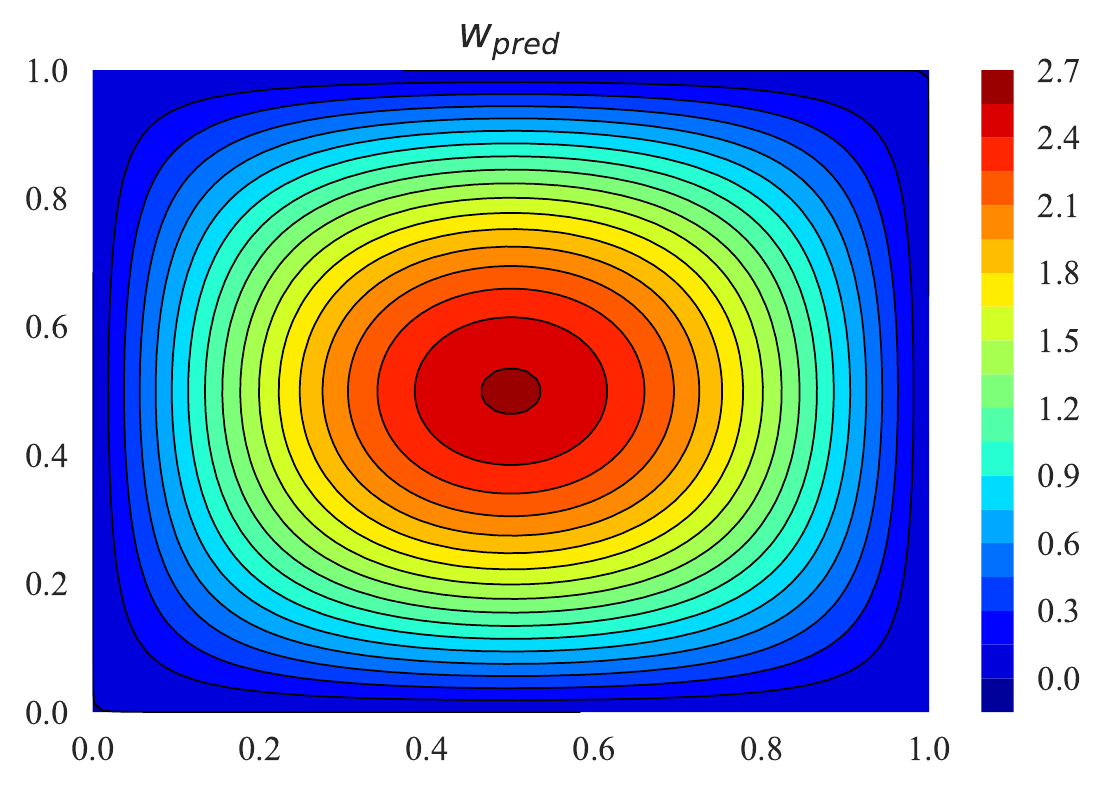}} & 
\hspace{-0.75cm}
\subfloat[Deflection error contour]{\includegraphics[width=7cm,height=6cm]{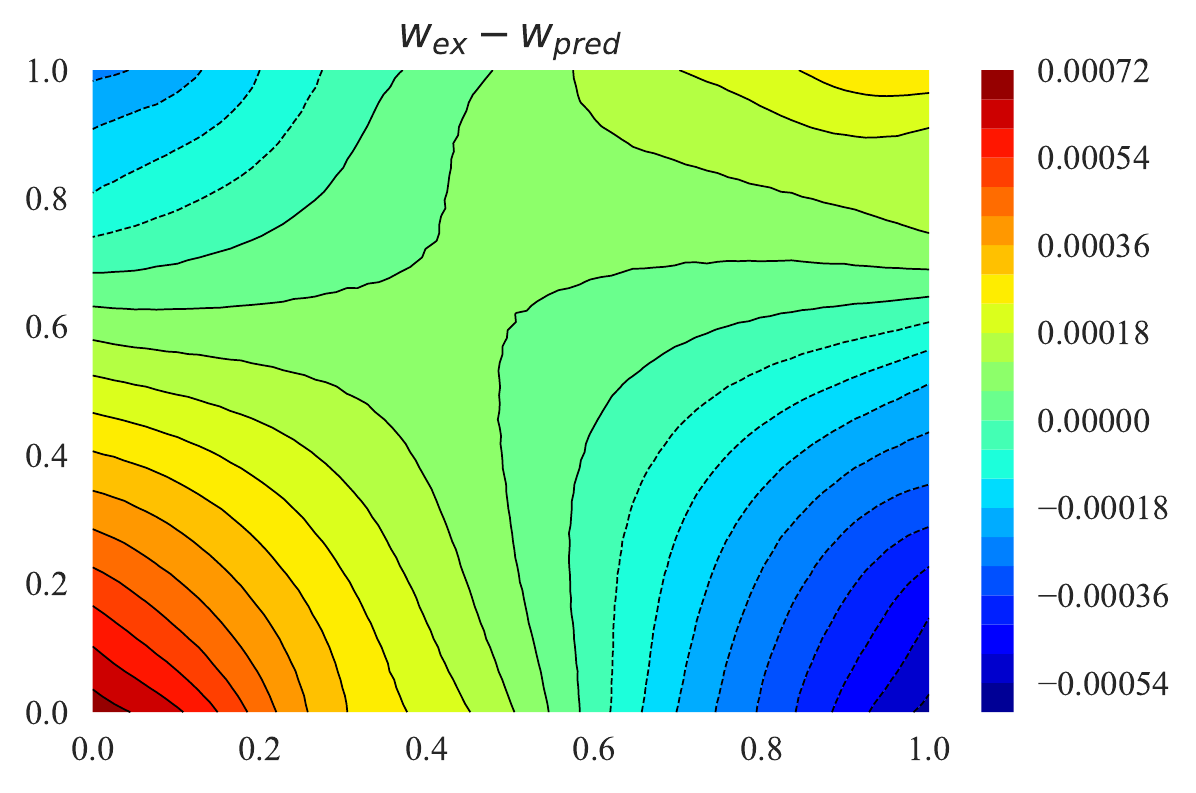}} \\ 
\vspace{-3cm}
\subfloat[Predicted deflection]{\includegraphics[width=7cm,height=6cm]{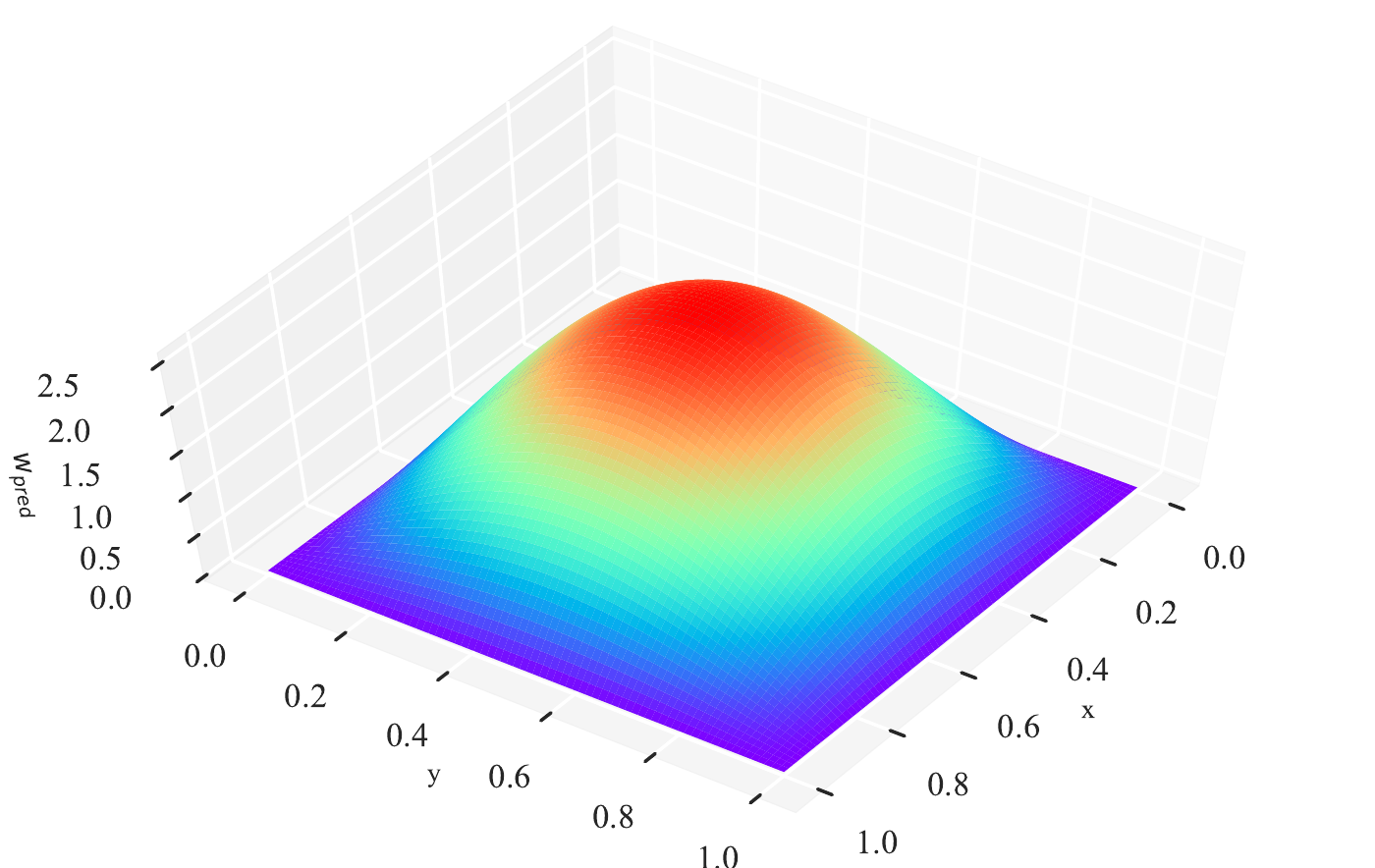}}&
\hspace{-0.75cm}
\subfloat[Exact deflection]{\includegraphics[width=7cm,height=6cm]{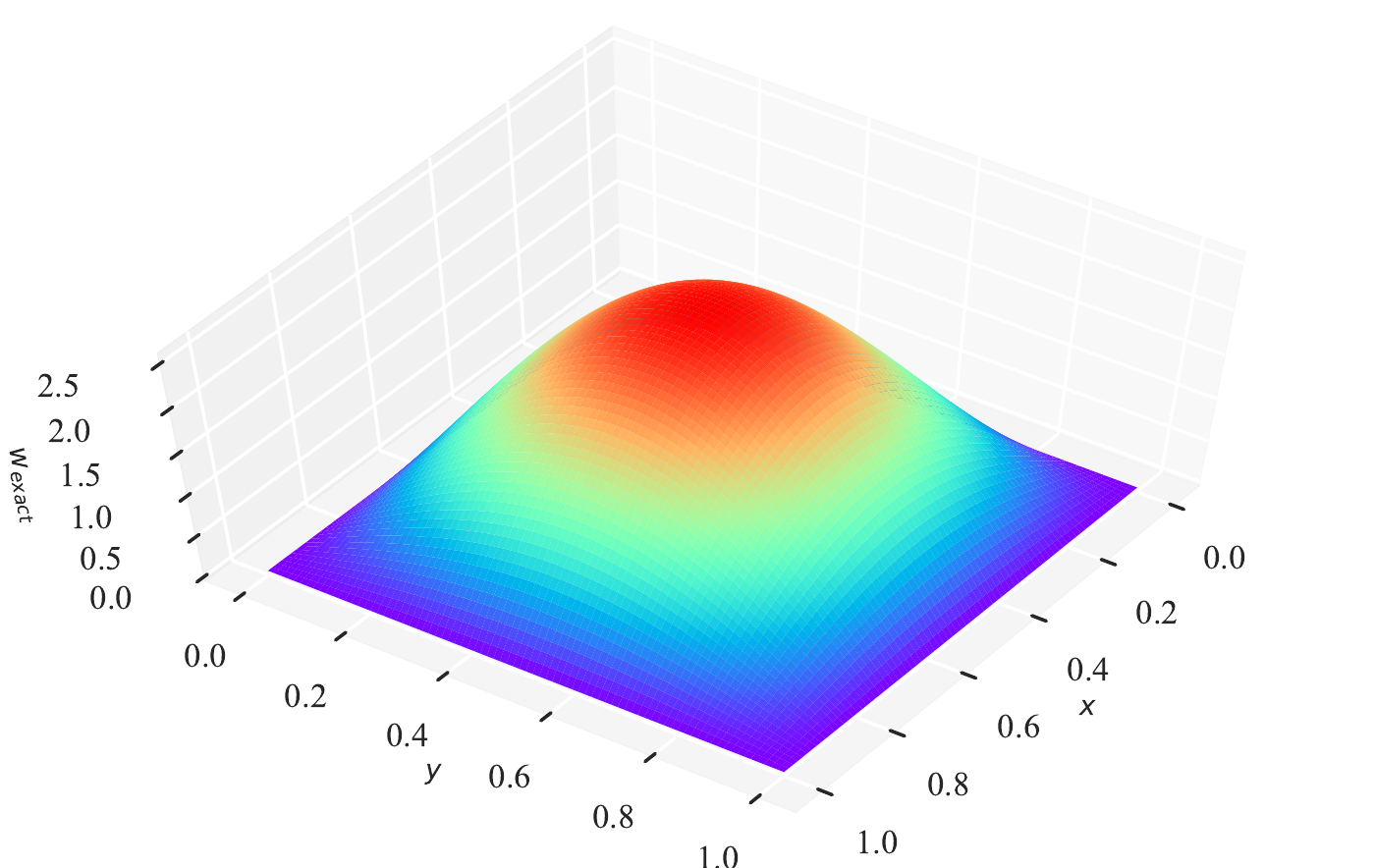}} \\
\end{tabular}
\vspace{3cm}
\caption{$\left(a\right)$ Predicted deflection contour $\left(b\right)$ Deflection error contour $\left(c \right)$ Predicted deflection $\left(d \right)$ Exact deflection of the simply-supported square plate with 2 hidden layers and 50 neurons.}
\label{Figure6:ssneuron50point2}
\end{figure}

\begin{figure}[H]
\centering
\begin{tabular}{cc}
\subfloat[Predicted deflection contour]{\includegraphics[width=7cm,height=6cm]{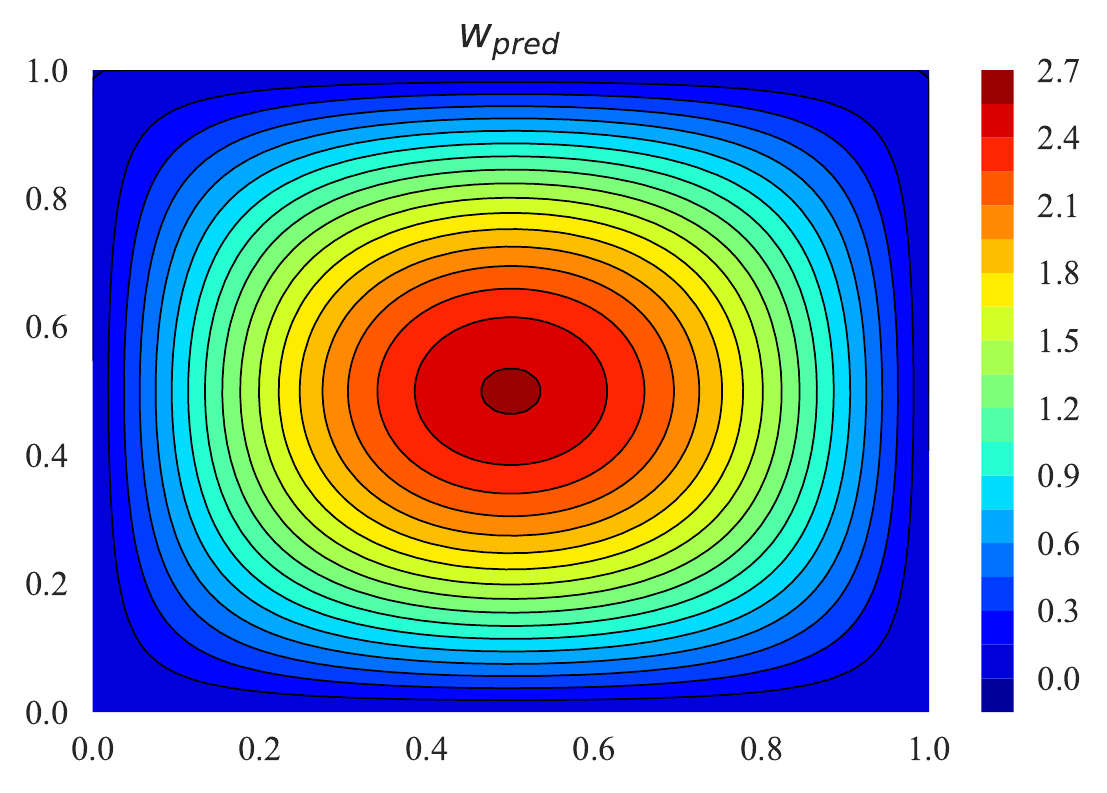}} & 
\hspace{-0.75cm}
\subfloat[Deflection error contour]{\includegraphics[width=7cm,height=6cm]{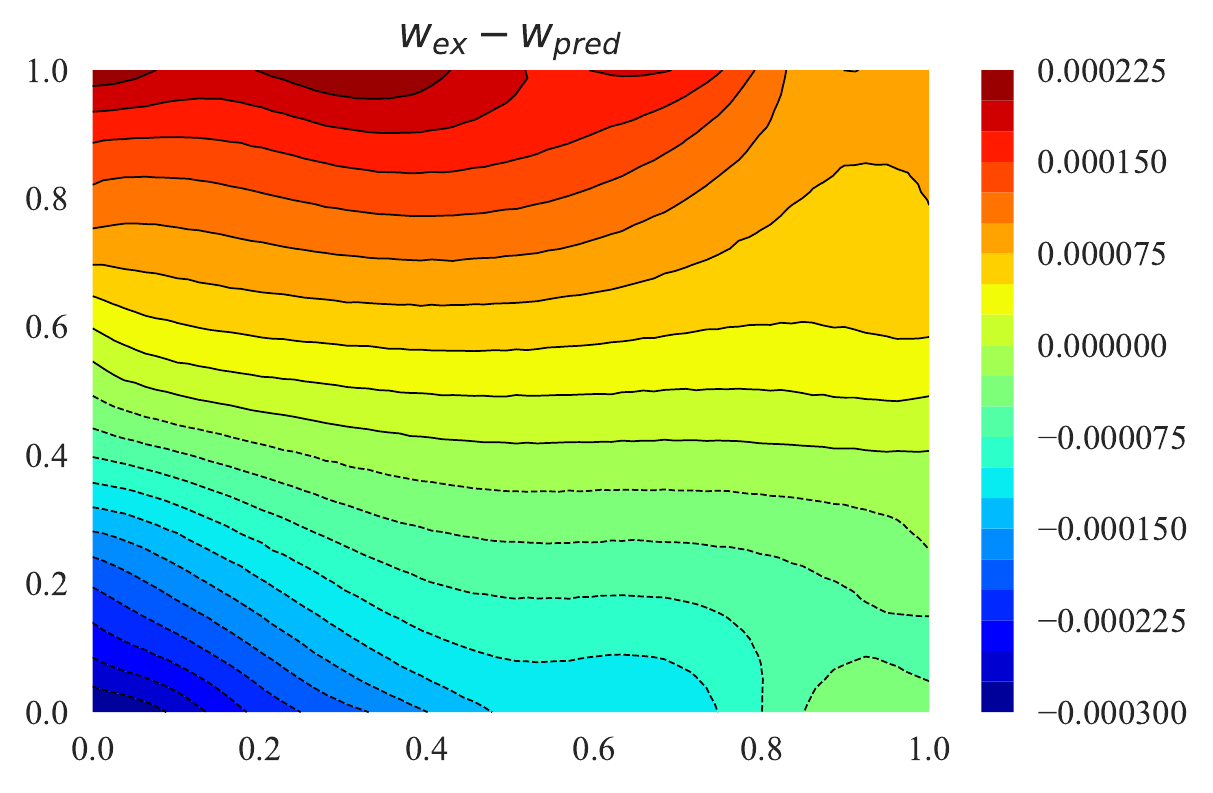}} \\ 
\vspace{-3cm}
\subfloat[Predicted deflection]{\includegraphics[width=7cm,height=6cm]{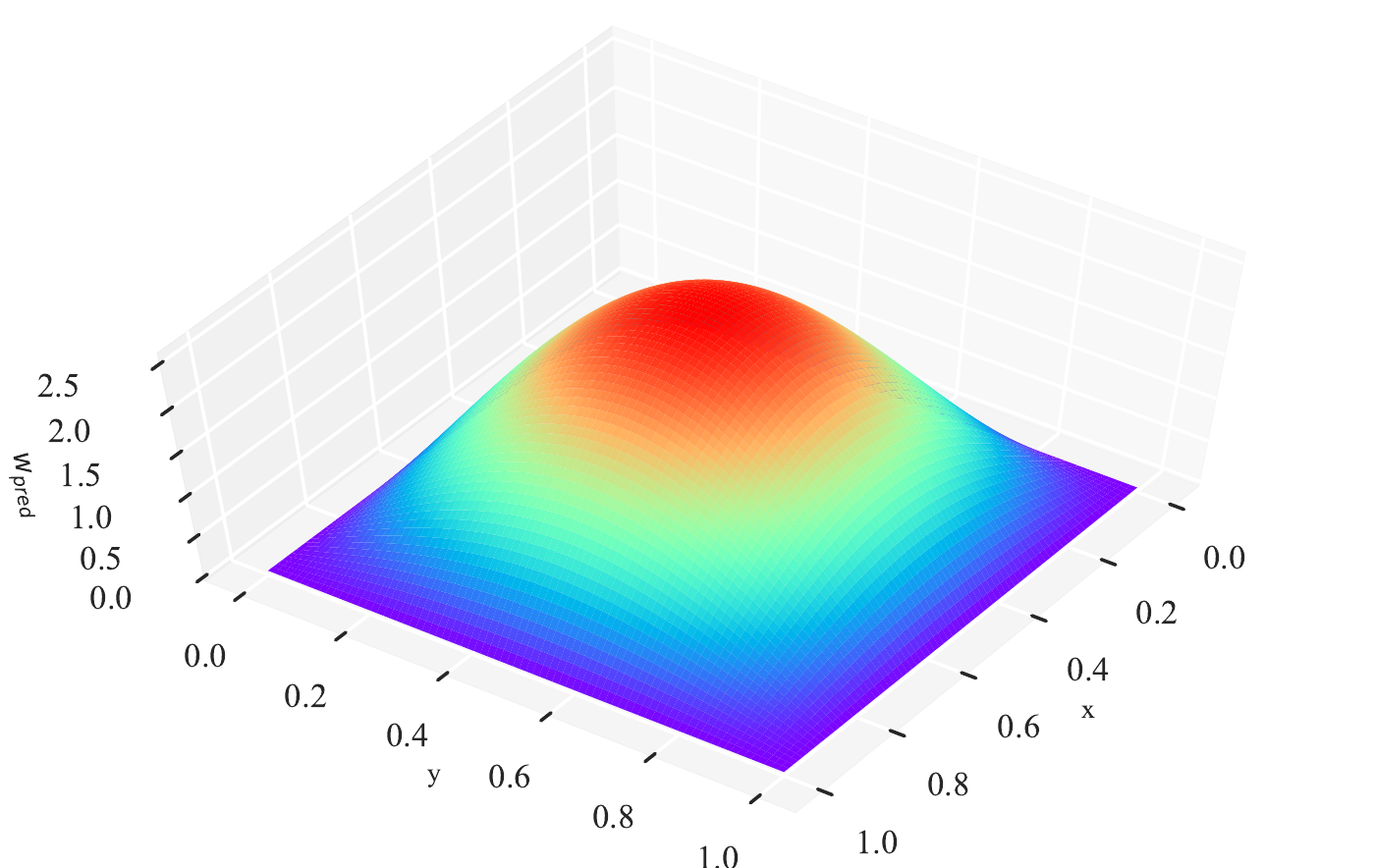}}&
\hspace{-0.75cm}
\subfloat[Exact deflection]{\includegraphics[width=7cm,height=6cm]{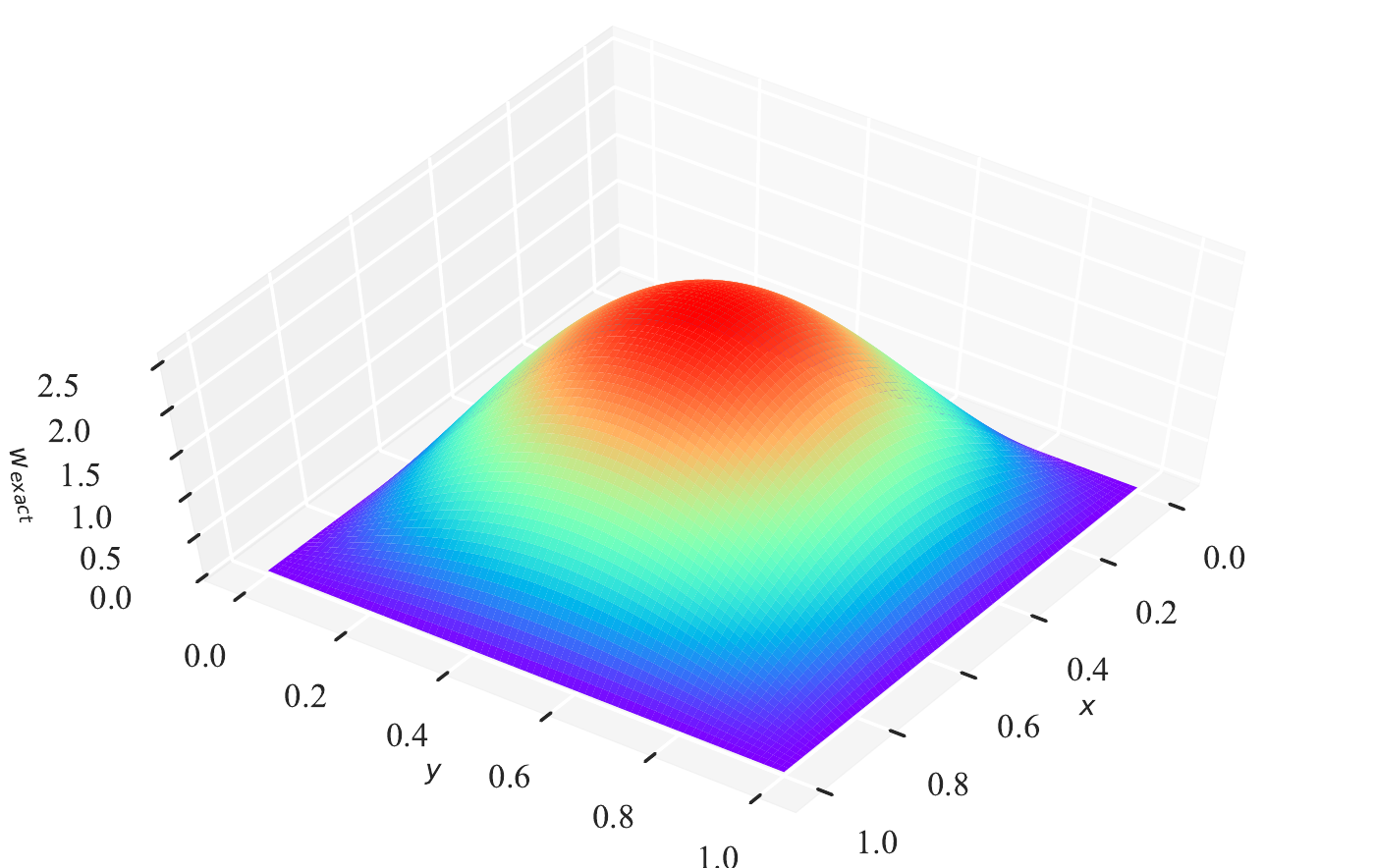}} \\
\end{tabular}
\vspace{3cm}
\caption{$\left(a\right)$ Predicted deflection contour $\left(b\right)$ Deflection error contour $\left(c \right)$ Predicted deflection $\left(d \right)$ Exact deflection of the simply-supported square plate with 3 hidden layers and 50 neurons.}
\label{Figure7:ssneuron50point3}
\end{figure}

\subsection{Clamped square plate}
\label{section 2:Clamped thin plate} 
A clamped square plate under a uniformly distributed transverse loading is also analyzed with deep collocation method in this section. There is no available explicit form exact solution for deflection of among the whole plate. And to better illustration the accuracy of this method, the analytical solution obtained by the Galerkin method referred in \cite{KHAN2012862} is adopted as a comparison:
\begin{equation}
\begin{bmatrix}
a_{11}\\ 
a_{12}\\ 
a_{21}\\ 
a_{22}\\ 
\end{bmatrix}=\frac{b^{4}p}{D}\begin{bmatrix}
0.318682766\\ 
0.038459815\\ 
0.038459815\\ 
0.008281438\\ 
\end{bmatrix},
\end{equation} 
\begin{equation}
\begin{matrix}
w=\frac{b^{4}q}{D}\left \{ a_{11}\left ( 1-\frac{x}{a} \right )^{2} \left ( 1-\frac{y}{b} \right )^{2} \left (\frac{x}{a} \right )^{2} \left (\frac{y}{b} \right )^{2}+
a_{12}\left ( 1-\frac{x}{a} \right )^{2} \left ( \frac{y}{b}-\frac{y^{2}}{b^{2}} \right )^{2} \left (\frac{x}{a} \right )^{2} \left (\frac{y}{b} \right )^{2}\right \}+\\ 
\frac{b^{4}q}{D}\left \{a_{21}\left ( \frac{x}{a} -\frac{x^{2}}{a^{2}}\right )^{2} \left ( 1-\frac{y}{b} \right )^{2} \left (\frac{x}{a} \right )^{2} \left (\frac{y}{b} \right )^{2}+
a_{22}\left ( \frac{x}{a} -\frac{x^{2}}{a^{2}}\right )^{2} \left (\frac{y}{b} -\frac{y^{2}}{b^{2}}\right )^{2} \left (\frac{x}{a} \right )^{2} \left (\frac{y}{b} \right )^{2}\right \}
\end{matrix}.
\end{equation}  

For the maximum transversal deflection at the centre of an isotropic square plate, Ritz method gives the maximum deflection at the centre as $w_{max}=0.00133\frac{qa^{4}}{D}$ \cite{KHAN2012862},and Timoshenko and Krieger \cite{timoshenko1959theory} gave a exact solution $w_{max}=0.00126\frac{qa^{4}}{D}$.  
Here, $D$ denotes the flexural stiffness of the plate
   and depends on the plate thickness and material properties. $a$,$b$ the length dimension of the plate. 1000 randomly generated collocation points as in Figure~\ref{Figure3:Scatterpoint} are also used to discritize the clamped square plate here. 

For this clamped case, a deep feedforward neural network with increasing layers and neurons is also studied in order to validate the convergence of this scheme. First, the maximum central deflection shown in Table 2 is also calculated for changing layers and neurons and are compared with aforementioned Ritz method, Galerkin method and exact solution by Timoshenko. It is demonstrated that our deep collocation method give most agreeable results with the exact solution. However, for neural networks with single hidden layer, the results are not that accurate even with 60 neurons. But as the neuron number increases, the results are indeed more accurate for the neural network with single hidden layer. This can be observed for the other two hidden layer types. Additionally, as the number of hidden layer increases, the results are much more accurate than the single hidden layer neural network results. 

The relative error with the analytical solution with different hidden layers and different neurons is shown in Figure~\ref{Figure8:Relativeerrorsc}. Although the magnitude of relative error of deflection for this numerical example is $1\times10^{-4}$, this dose not mean that our deep collocation method is not that accurate for this problem. For it is mentioned that the deflection vector as a comparison to calculate the relative error is gained from Galerkin method, and we have gotten from Table 2, that our method gives more accurate maximum deflection than Galerkin method. As hidden layers increase, the two flat relative error curves nearly coincide and converge to the exact solution.
\begin{figure}[H]
\centering
\caption*{Table 2: Maximum deflection predicted by deep collocation method.}
\vspace{-0.55cm}
\begin{tabular}{c}
\subfloat{\includegraphics[width=0.95\textwidth]{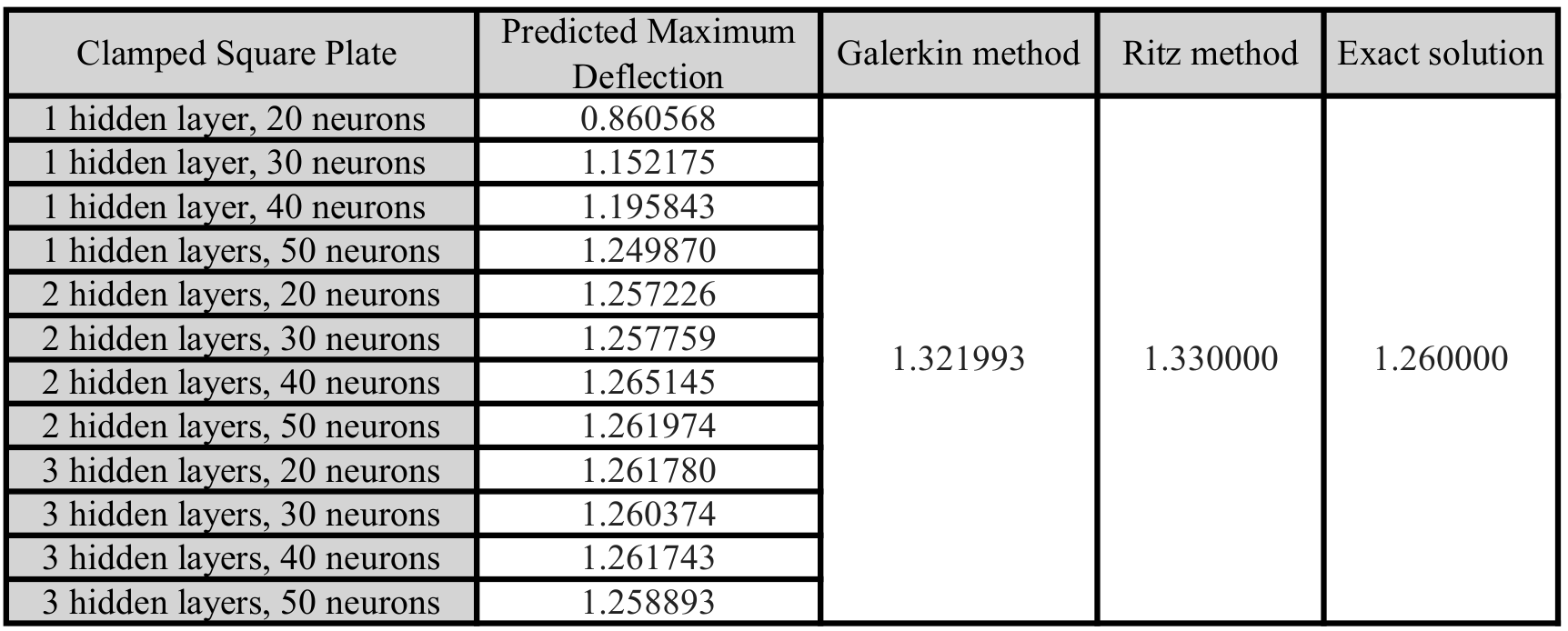}}
\end{tabular}
\end{figure}

\begin{figure}[H]
\centering
\begin{tabular}{c}
\subfloat{\includegraphics{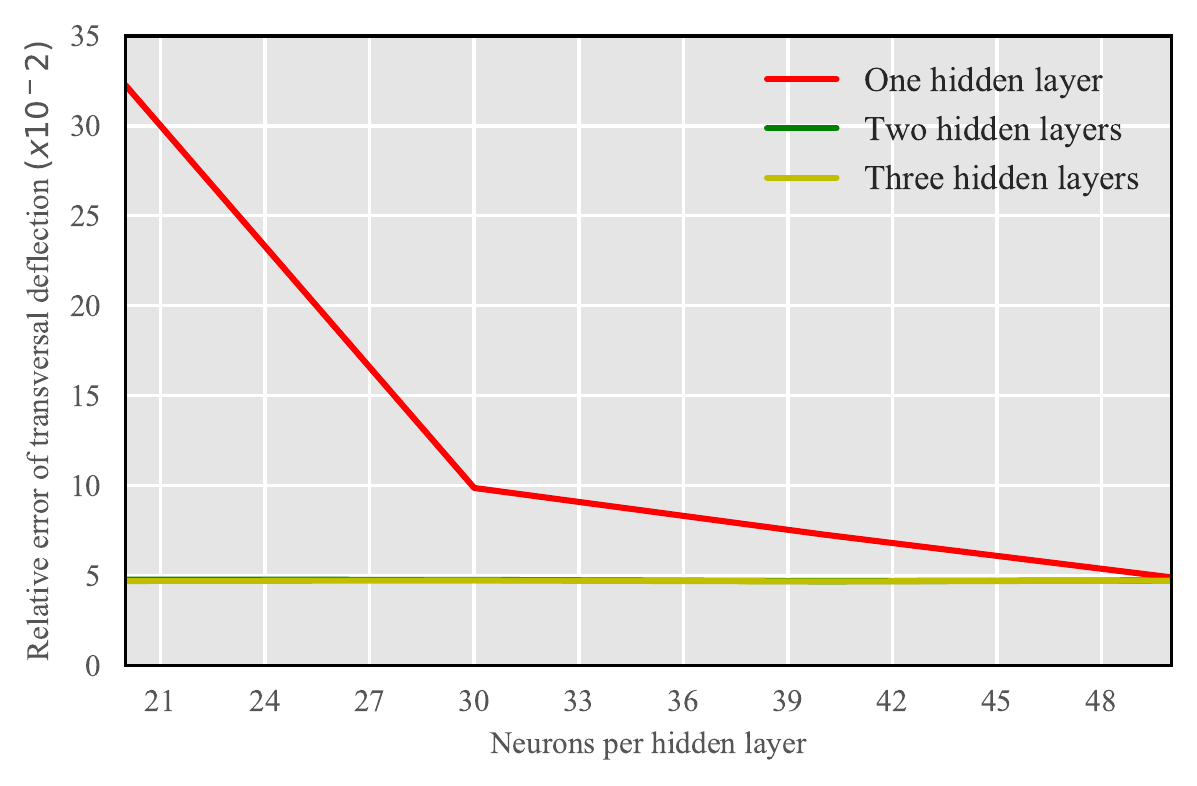}}
\end{tabular}
\vspace{-0.55cm}
\caption{The relative error of deflection with varying hidden layers and neurons.}
\label{Figure8:Relativeerrorsc}
\end{figure}

\begin{figure}[H]
\centering
\begin{tabular}{cc}
\subfloat[Predicted deflection contour]{\includegraphics[width=7cm,height=6cm]{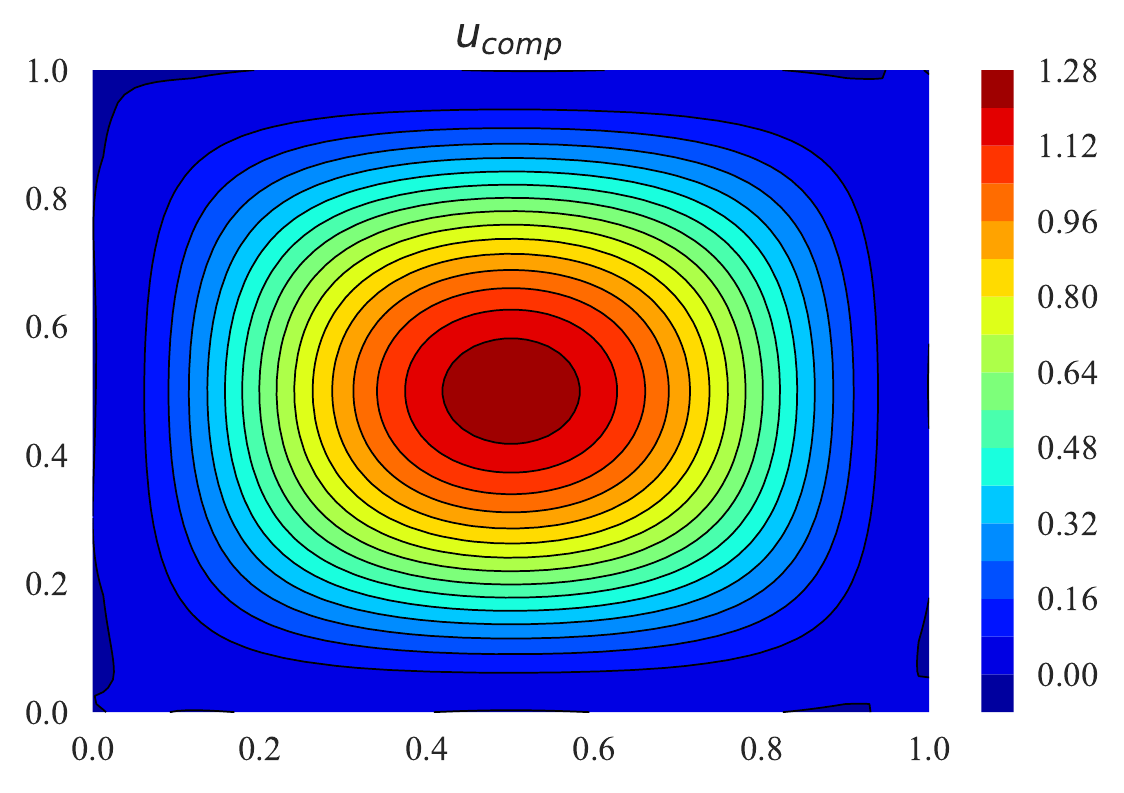}} & 
\hspace{-0.75cm}
\subfloat[Deflection error contour]{\includegraphics[width=7cm,height=6cm]{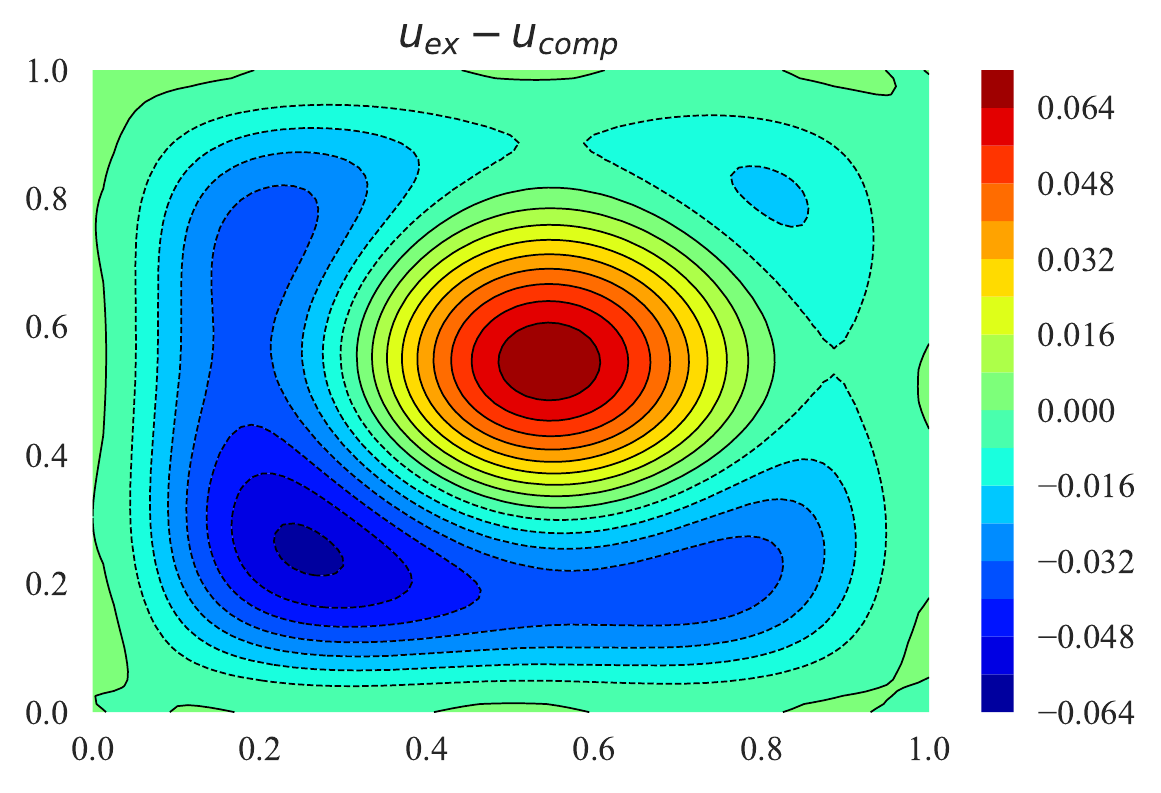}} \\ 
\vspace{-3cm}
\subfloat[Predicted deflection]{\includegraphics[width=7cm,height=6cm]{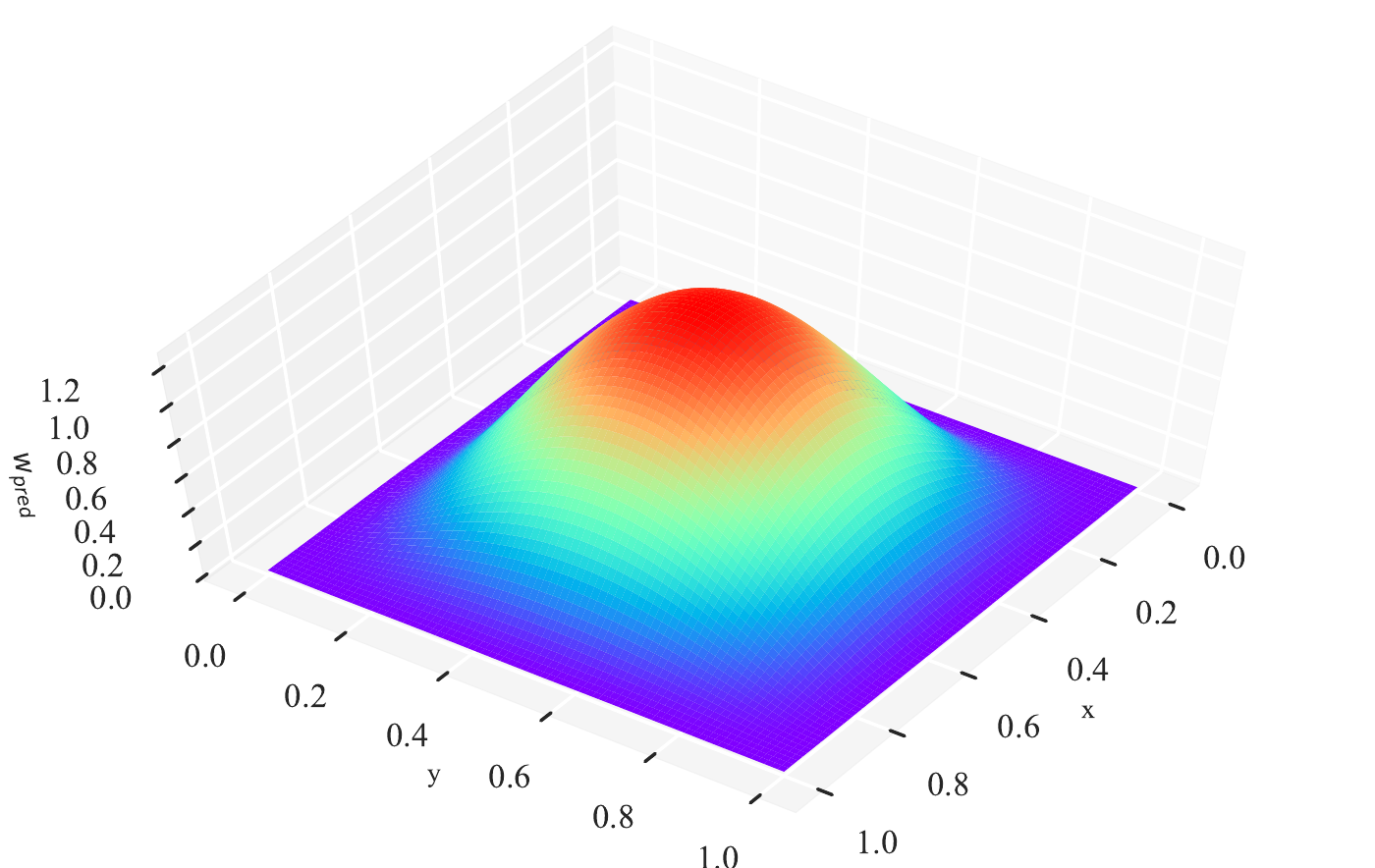}}&
\hspace{-0.75cm}
\subfloat[Exact deflection]{\includegraphics[width=7cm,height=6cm]{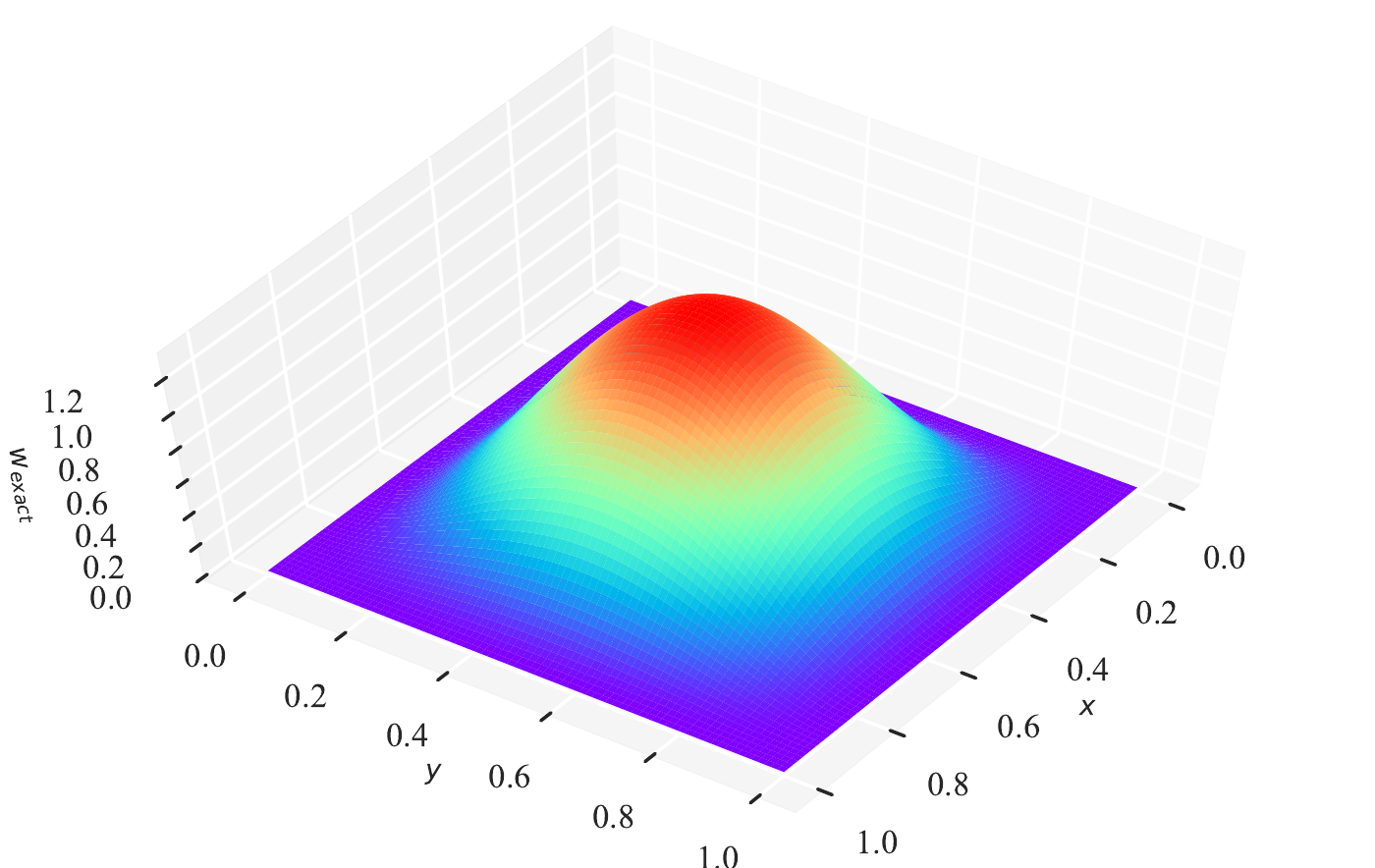}} \\
\end{tabular}
\vspace{3cm}
\caption{$\left(a\right)$ Predicted deflection contour $\left(b\right)$ Deflection error contour $\left(c \right)$ Predicted deflection $\left(d \right)$ Exact deflection of the clamped square plate with 3 hidden layers and 50 neurons.}
\label{Figure9:csneuron50point1}
\end{figure}

Finally, to better depict the favourable of our method, the deflection contour, relative error contour and deformed deflection of the middle surface are also listed for the deep neural network with three layers and 50 neurons in Figure~\ref{Figure9:csneuron50point1}. It is clear that the deep collocation method yields results agrees well with the analytical solution.

\subsection{Clamped circular plate}
\label{section 1:Circular thin plate} 
\begin{figure}[H]
\centering
\begin{tabular}{c}
\subfloat{\includegraphics{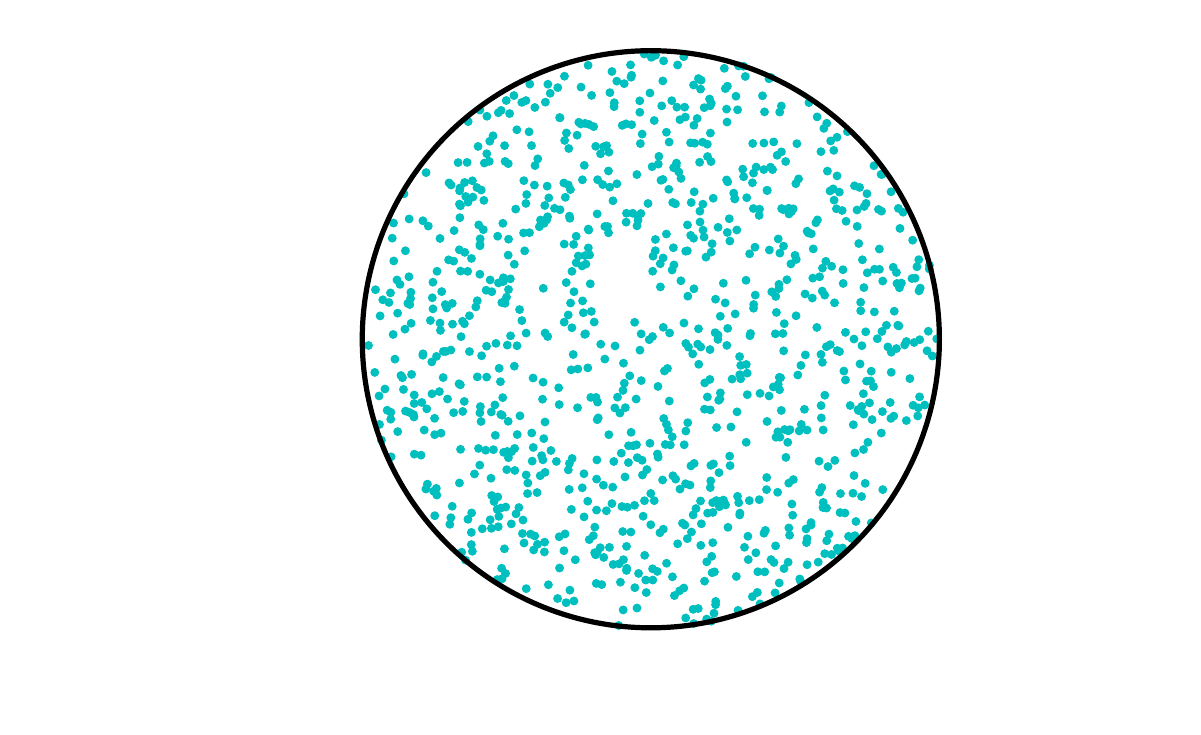}}
\end{tabular}
\vspace{-1.25cm}
\caption{Collocation points discretize the circular domain.}
\label{Figure10:Scatterpointcc}
\end{figure}

A clamped circular plate with radius $R$ under a uniform load $p$ is studied here. 1000 collocation points shown in Figure~\ref{Figure10:Scatterpointcc} are deployed among the circular plate first. Then, we applied deep collocation method to study the deformation of this circular plate. This problem has a exact solution, which can be referred in \cite{timoshenko1959theory}:
\begin{equation}
w=\frac{p\left ( R^{2}-\left ( x^{2}+ y^{2}\right ) \right )^{2}}{64D},
\end{equation}   
with $D$ denotes the flexural stiffness of the plate
   and depends on the plate thickness and material properties.

The maximum deflection at the central of the circular plate with varying hidden layers and neurons in Table 3 and compared with exact solution. It is obvious that the predicted maximum deflection is very accurate, and as the neuron and hidden number increase, the maximum deflection are more and more close to the exact solution.

The relative error for deflection of clamped circular plate with increasing hidden layers and neurons is shown in Figure~\ref{Figure11:Relativeerrorcc} in order to show the convergent of this method. From this figure, we can get that the as hidden layer number increases, the relative error curves become flat and converge very well to the exact solution. However, all neural networks perform well with a relative error magnitude of $1\times10^{-4}$.

\begin{figure}[H]
\centering
\caption*{Table 3: Maximum deflection predicted by deep collocation method.}
\vspace{-0.55cm}
\begin{tabular}{c}
\subfloat{\includegraphics{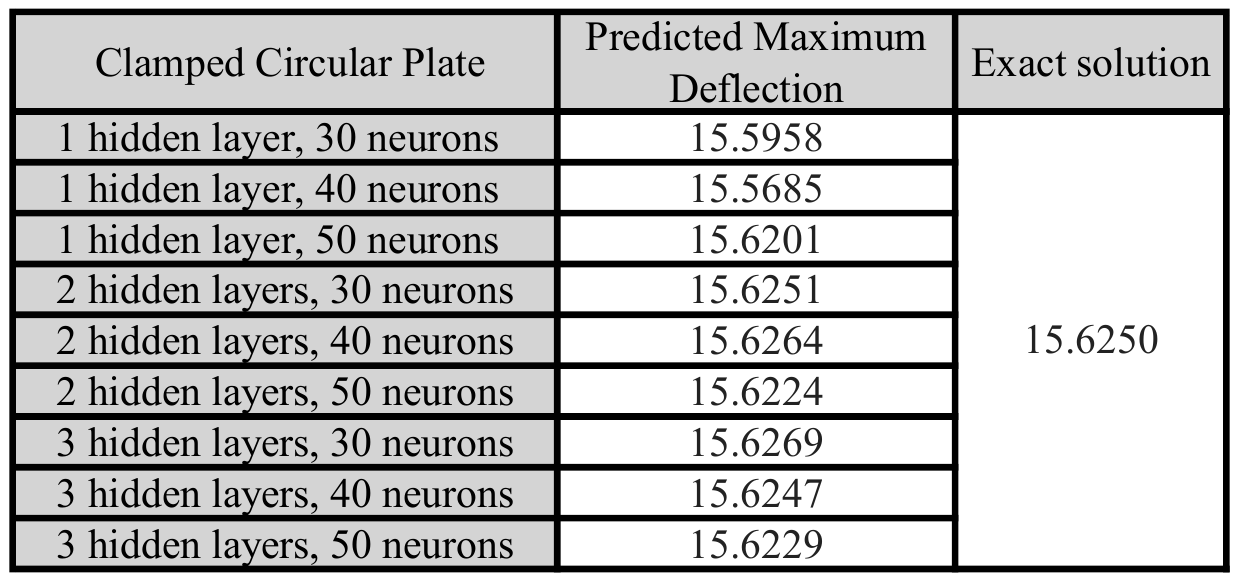}}
\end{tabular}
\end{figure}

Finally, the deformation contour, deflection error contour, predicted and exact deformation figure are displayed in Figure~\ref{Figure12:ccneuron50point1}. The deflection of this circular plate agrees well with the exact solution. The accuracy of this collocation method is again shown here, which also illustrates that this deep collocation method can be easily and agreeably applied to simulate deformation of plates of various shapes.

\begin{figure}[H]
\centering
\begin{tabular}{c}
\subfloat{\includegraphics{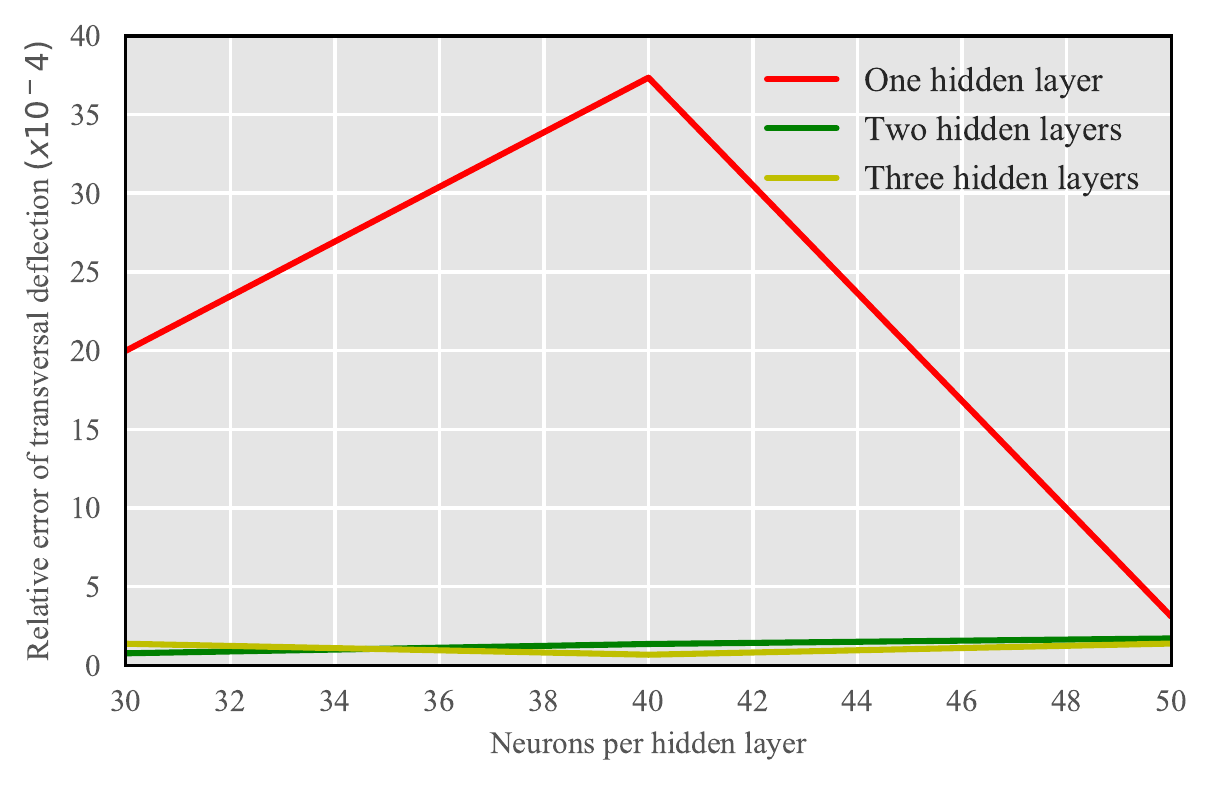}}
\end{tabular}
\vspace{-0.55cm}
\caption{The relative error of deflection with varying hidden layers and neurons.}
\label{Figure11:Relativeerrorcc}
\end{figure}

\begin{figure}[H]
\centering
\begin{tabular}{cc}
\subfloat[Predicted deflection contour]{\includegraphics[width=7cm,height=6cm]{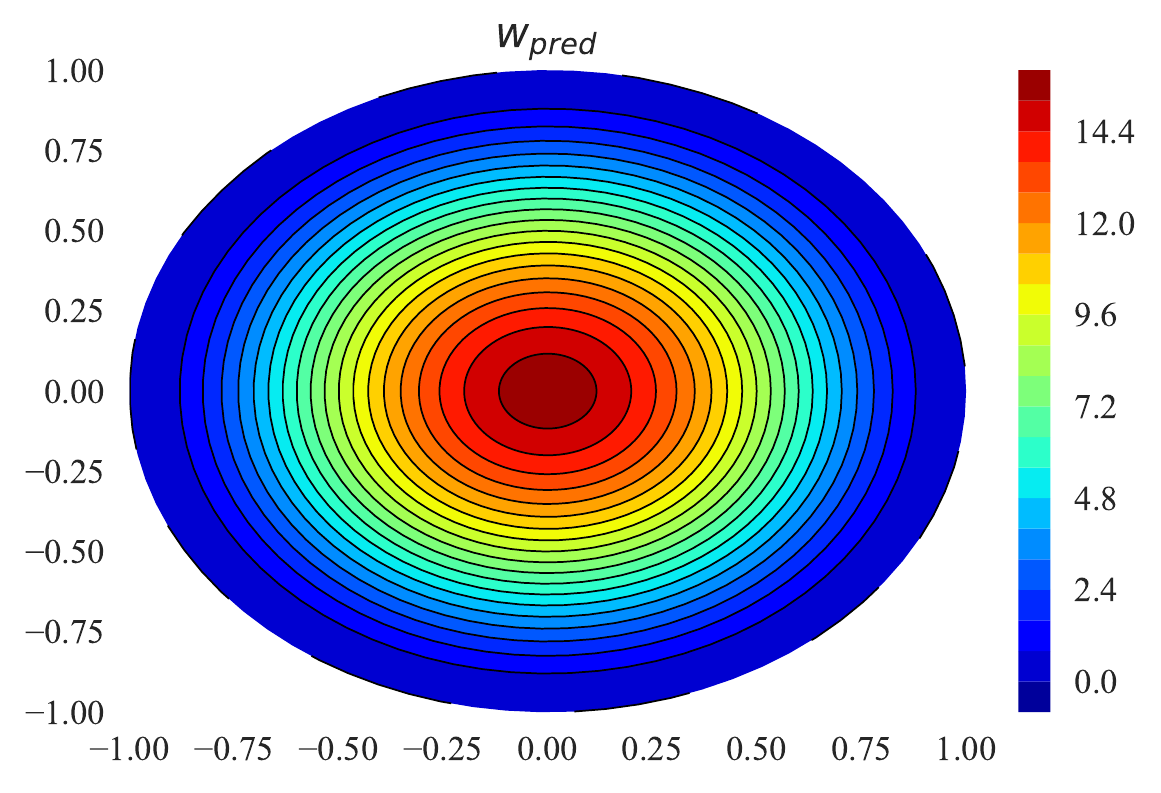}} & 
\hspace{-0.75cm}
\subfloat[Deflection error contour]{\includegraphics[width=7cm,height=6cm]{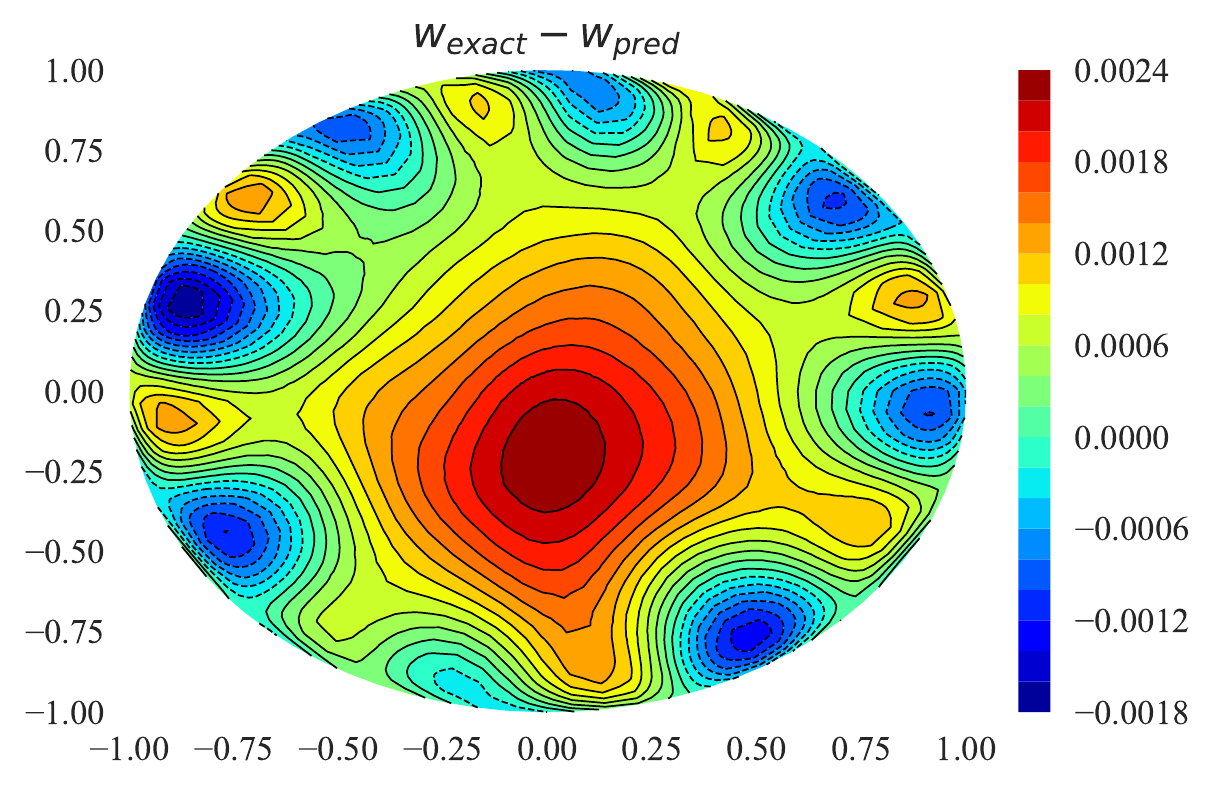}} \\ 
\vspace{-3cm}
\subfloat[Predicted deflection]{\includegraphics[width=7cm,height=6cm]{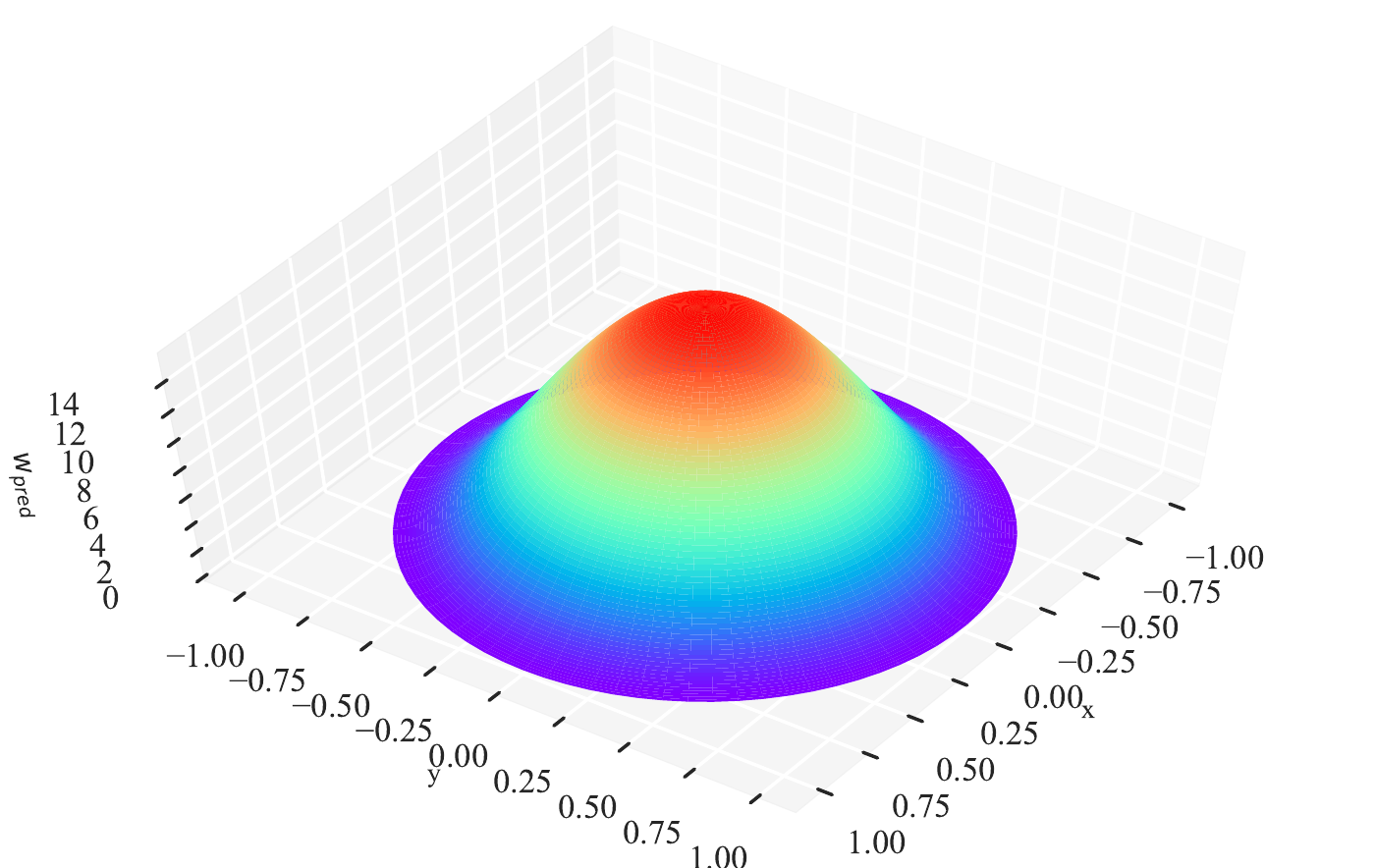}}&
\hspace{-0.75cm}
\subfloat[Exact deflection]{\includegraphics[width=7cm,height=6cm]{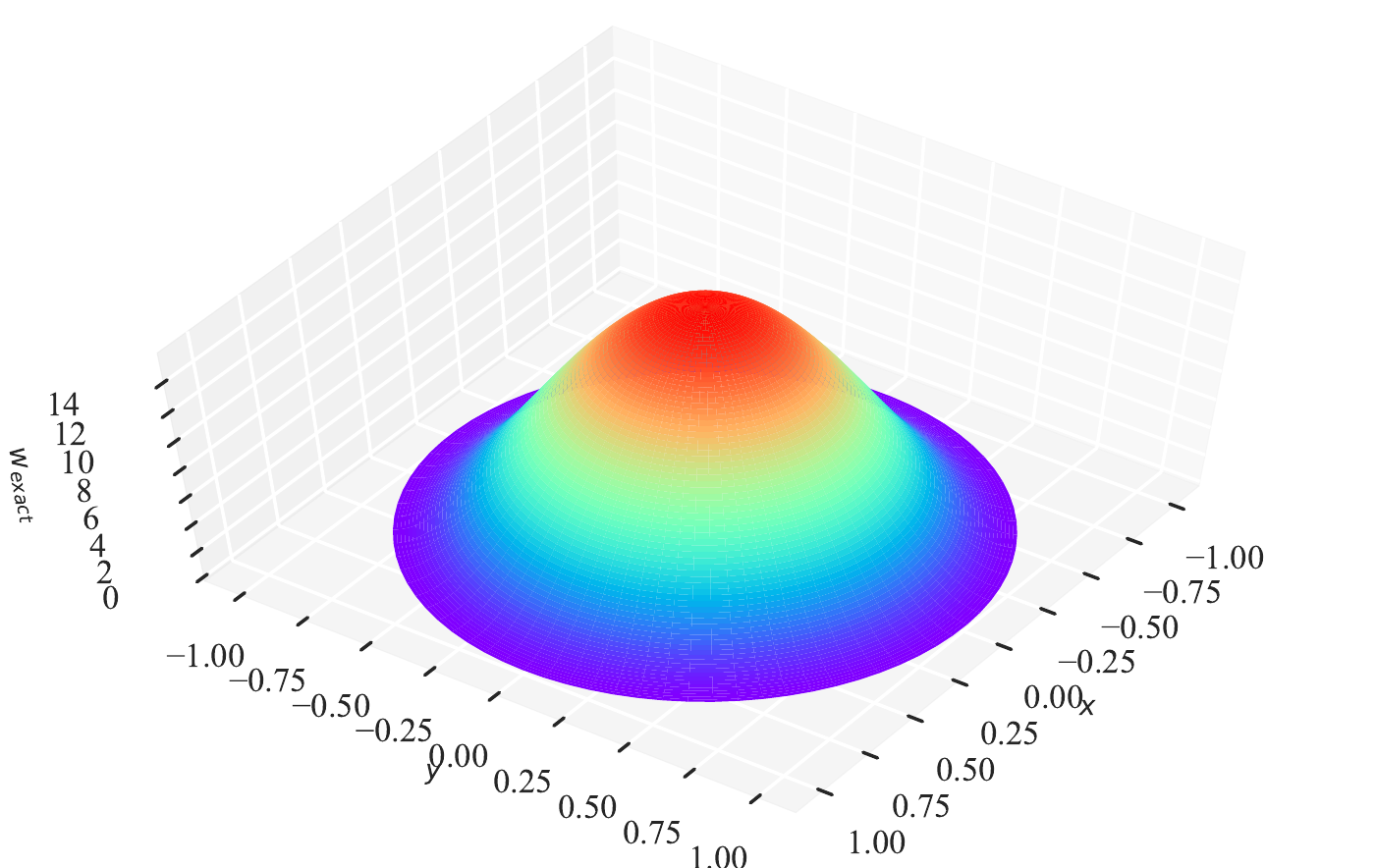}} \\
\end{tabular}
\vspace{3cm}
\caption{$\left(a\right)$ Predicted deflection contour $\left(b\right)$ Deflection error contour $\left(c \right)$ Predicted deflection $\left(d \right)$ Exact deflection of the clamped circular plate with 3 hidden layers and 50 neurons.}
\label{Figure12:ccneuron50point1}
\end{figure}

\subsection{Simply-supported square plate on Winkler foundation}
\label{section 1:simply-supported square plate} 
The simply-supported square plate resting on Winkler foundation is studied in this section, which assumes that the foundation's reaction $p\left(x,y\right)$ can be described by $p\left(x,y\right)=\textit{k}w$, with $\textit{k}$ a constant called $foundation\;modulus$. Considering a plate on a continuous Winkler foundation, the governing Equation \ref{governing} can be written as
\begin{equation}
\bigtriangledown^{2}\left ( \bigtriangledown^{2}w \right )=\bigtriangledown^{4}w=\frac{\left(p-q\right)}{D}=\frac{\left(p-\textit{k}w\right)}{D}
\label{governingelastfound}
\end{equation}   
The analytical solution for this numerical example is given as \cite{timoshenko1959theory}:
\begin{equation}
w=\frac{16p}{ab}\sum_{m=1,3,5,\cdots  }^{\infty}\sum_{n=1,3,5,\cdots  }^{\infty}\frac{\textrm{sin}\frac{m\pi x }{a} \textrm{sin}\frac{n\pi y }{b}}{mn\left [ \pi^4D\left ( \frac{m^2}{a^2} +\frac{n^2}{b^2}\right )^2+\textit{k} \right ]}
\end{equation}   

For this numerical example, the arrangement of collocation points are the same as that in Figure~\ref{Figure3:Scatterpoint}. For the detail implementation, neural networks with different neurons and deepth are applied in the calculation. Also, maximum deflections shown in Table 4 at the central point in all those cases are first studied in order to unveil the accuracy of the deep collocation method. 

Good agreement can be obsevered in this numerical example as well. From Table 4, we can obsevered as hidden layer and neuron number grows, the maximum deflection becomes more accurate and close to the analytical serial solution for even two hidden layers. The relative error shown in Figure~\ref{Figure13:Relativeerrorssef}  better depicts the advantages of deep neural network than shallow wide neural network. And with more hidden layers, with neurons increase, the relative error cure becomes flat and very close to zero, which shows that the deep collocation method with only two hidden layers can well approximate the deflection.

To better illustrate the deflection distribution around the whole plate, deflection contour, deflection error contour, deformation contour on deformed figure are shown in Figure~\ref{Figure14:esneuron50point1} and compared with the analytical solution. It is demonstrated that the proposed method agrees well with the analytical solution.

\begin{figure}[H]
\centering
\caption*{Table 4: Maximum deflection predicted by deep collocation method.}
\vspace{-0.55cm}
\begin{tabular}{c}
\subfloat{\includegraphics{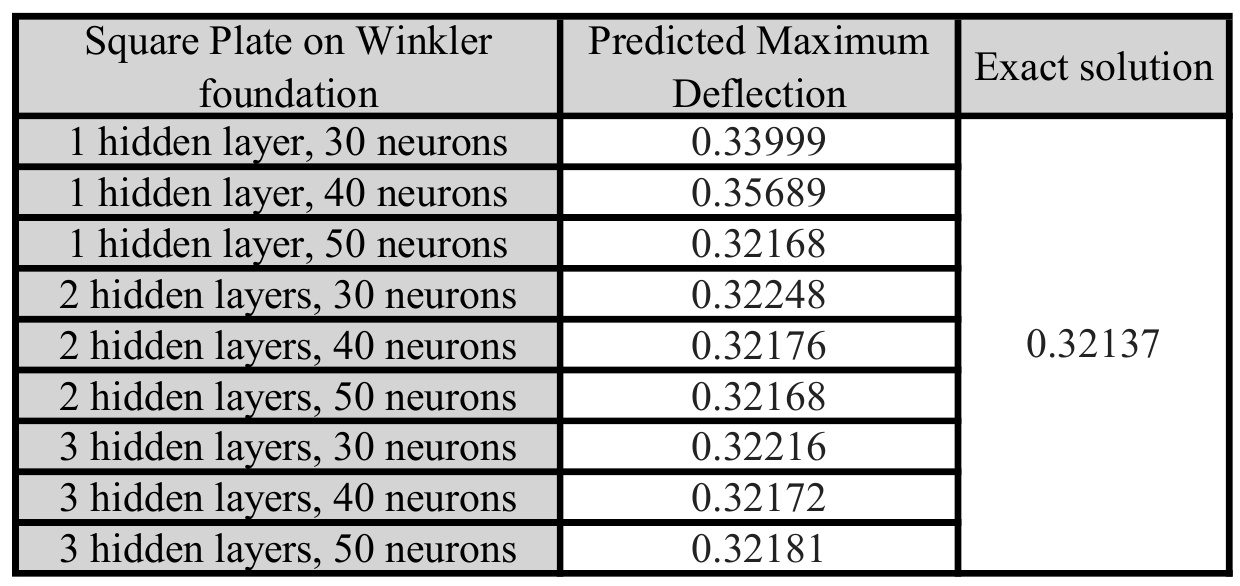}}
\end{tabular}
\end{figure}

\begin{figure}[H]
\centering
\begin{tabular}{c}
\subfloat{\includegraphics{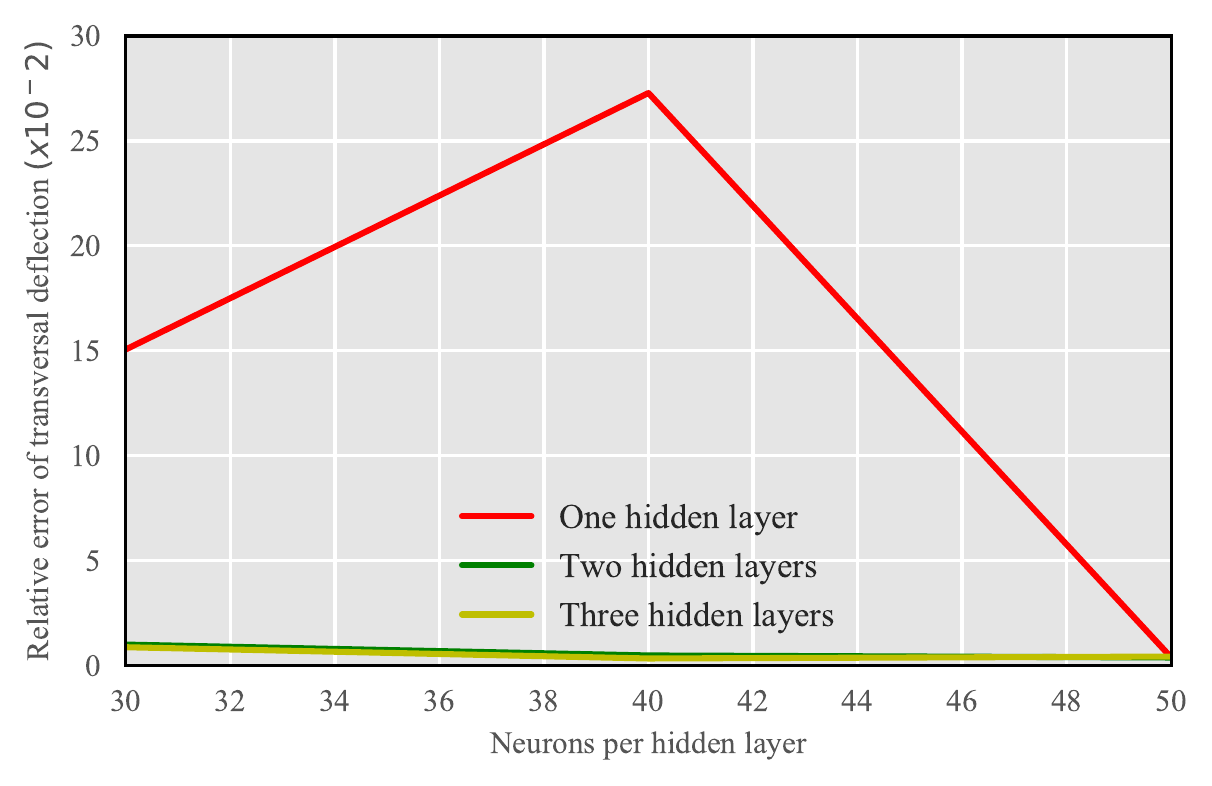}}
\end{tabular}
\vspace{-0.55cm}
\caption{The relative error of deflection with varying hidden layers and neurons.}
\label{Figure13:Relativeerrorssef}
\end{figure}

\begin{figure}[H]
\centering
\begin{tabular}{cc}
\subfloat[Predicted deflection contour]{\includegraphics[width=7cm,height=6cm]{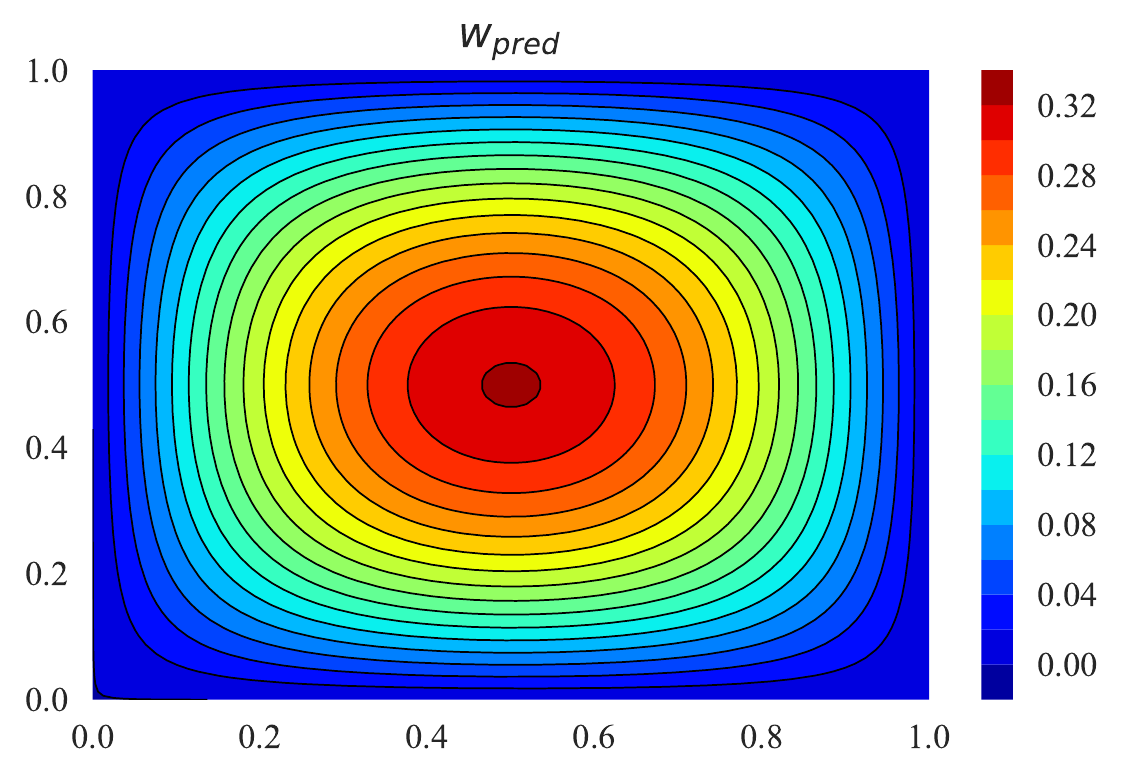}} & 
\hspace{-0.75cm}
\subfloat[Deflection error contour]{\includegraphics[width=7cm,height=6cm]{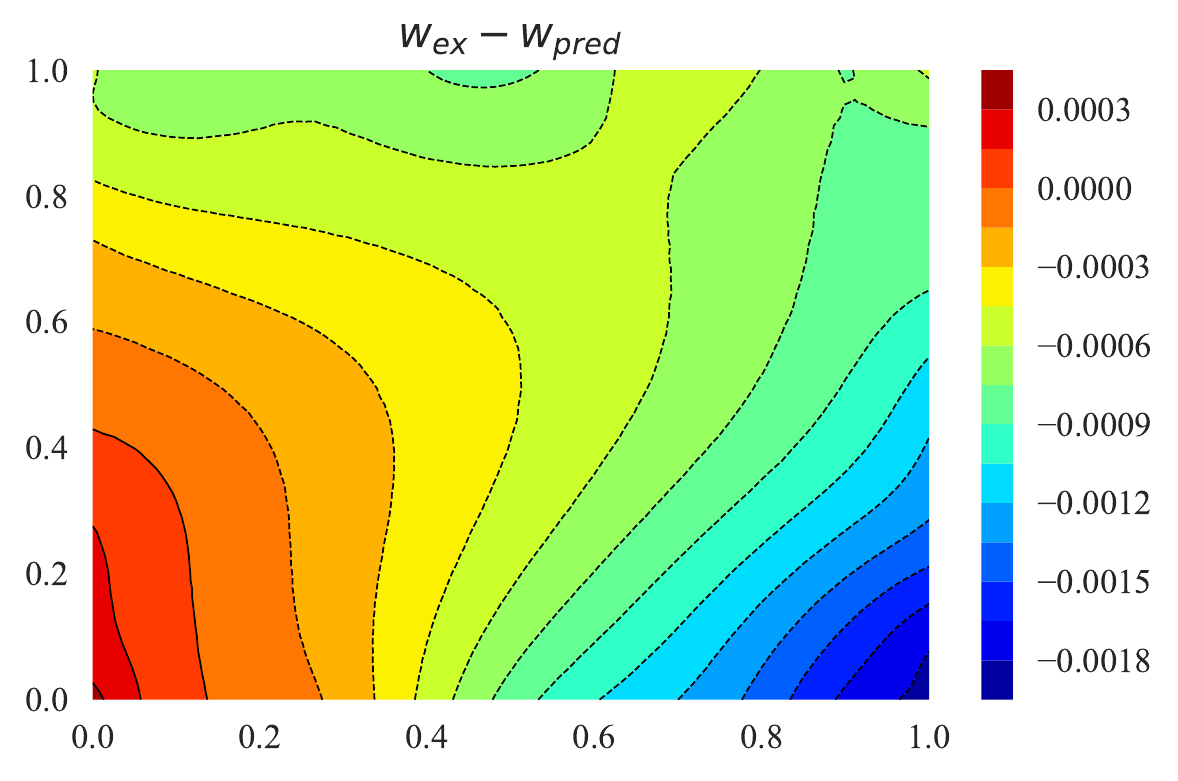}} \\ 
\vspace{-3cm}
\subfloat[Predicted deflection]{\includegraphics[width=7cm,height=6cm]{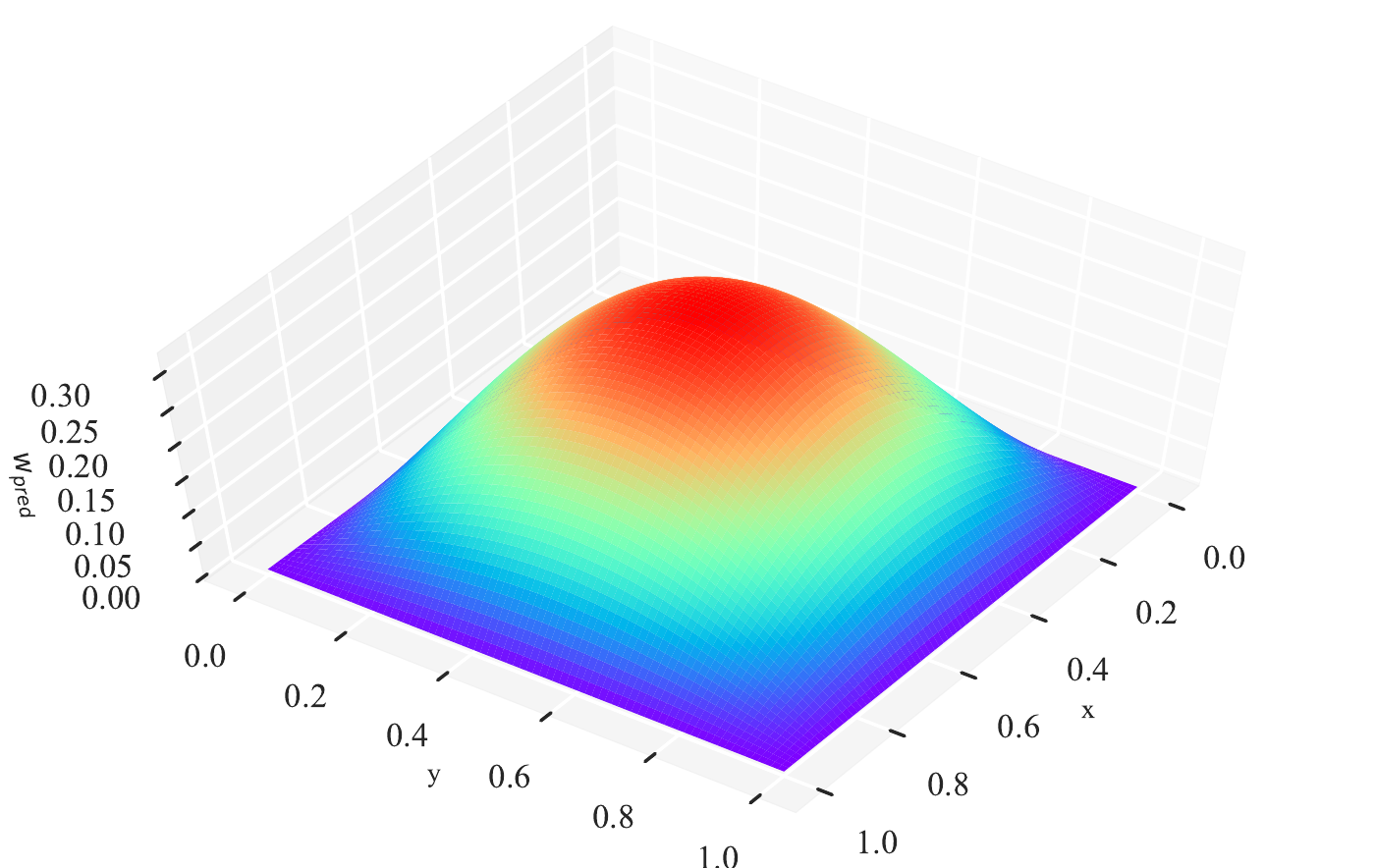}}&
\hspace{-0.75cm}
\subfloat[Exact deflection]{\includegraphics[width=7cm,height=6cm]{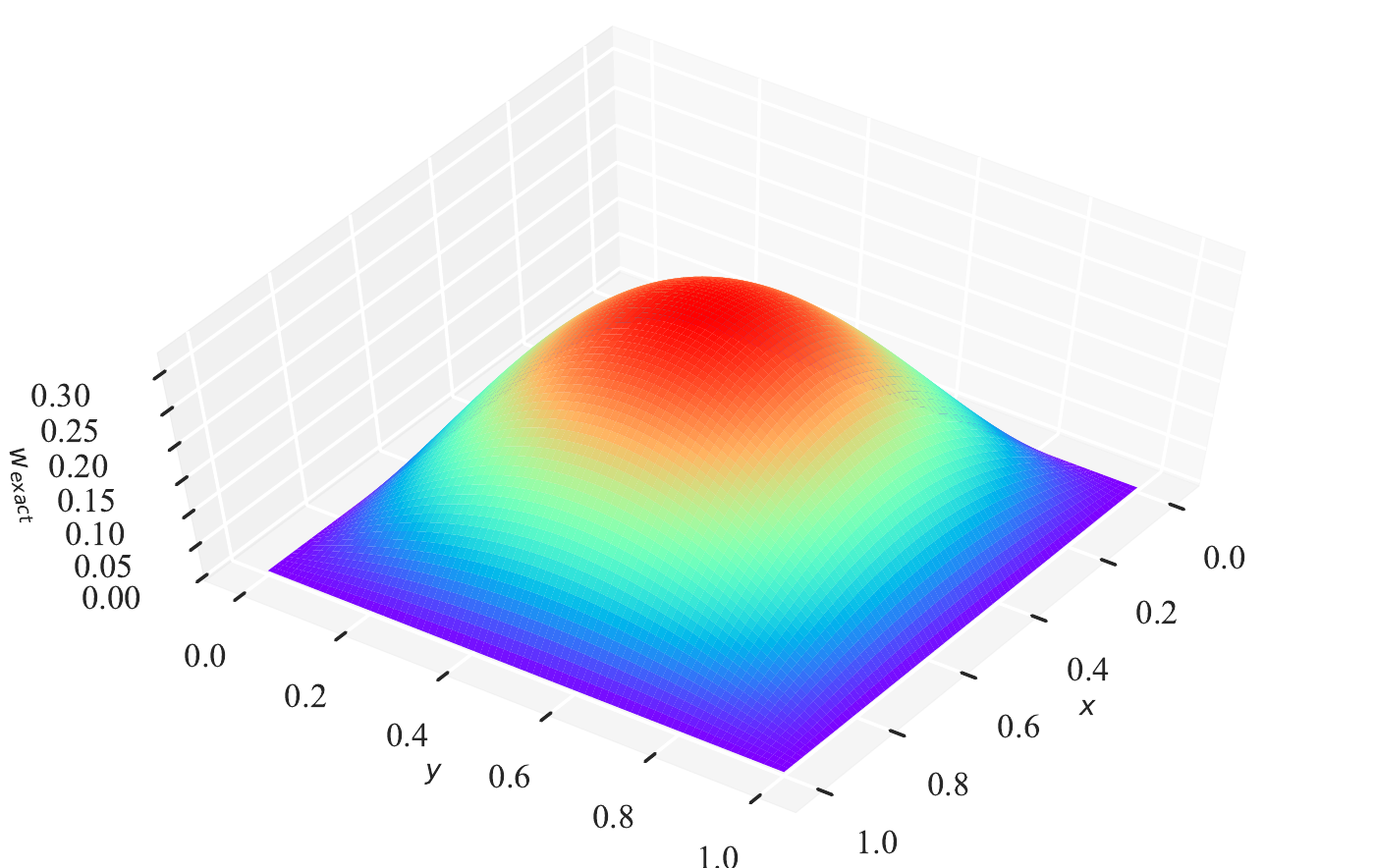}} \\
\end{tabular}
\vspace{3cm}
\caption{$\left(a\right)$ Predicted deflection contour $\left(b\right)$ Deflection error contour $\left(c \right)$ Predicted deflection $\left(d \right)$ Exact deflection of the simply-supported plate on Winkler foundation with 3 hidden layers and 50 neurons.}
\label{Figure14:esneuron50point1}
\end{figure}

\section{Conclusions}
In this study,study the bending analysis of Kirchhoff plates of various shapes, loads and boundary conditions. The governing equation of this problem is a fourth order partial differential equation (biharmonic equation), which is an important kind of PDEs in engineering mechanics. The proposed deep collocation method is a truly "meshfree" method, and can be used to approximate any continuous function, which is very suitable for the analysis of thin plate bending problems. The deep collocation method is very simple in implementation, which can be further applied in a wide variety of engineering problems.

Moreover, the deep collocation method with randomly distributed collocations and deep neural networks perform very well with a MSE loss function minimized by the combined L-BFGS and Adam optimizer. An accurate result can even be gotten for the single layer and 20 neurons case. However, as the increase of hidden layers and neurons on each layer, most results become more accurate and converge to the exact and analytical solution. For circular plates, this method become extremely efficient and accurate, and accurate results can be obtained with only a few layers and neurons. More importantly, once those deep neural networks are trained, they can be used to evaluate the solution at any desired points with minimal additional computation time.

However, there are still some intriguing issues remained to be studied for the deep neural network based method such as the influence of choosing other neural network types, activation functions, loss function forms, weight/bias initialization, and optimizers on the accuracy and efficiency of this deep collocation method, which will be studied in our future research.

\bigskip
\noindent\textbf{Acknoledgement}:

\bibliographystyle{unsrt}
\bibliography{ref/DLplate.bib}

\end{document}